\documentclass[12pt,leqno]{article}

\usepackage{amsfonts}
\usepackage{amsmath}
\usepackage{epsfig}
\usepackage{graphicx}
\input epsf
\usepackage{latexsym}
\usepackage[latin1]{inputenc}
\usepackage[OT1]{fontenc}

\usepackage{color}
\definecolor{purple}{rgb}{0.65, 0, 1}
\definecolor{orange}{rgb}{1,.5,0}
\definecolor{brown}{rgb}{.9,.73,.26}

\def\RR{\mathbb R}
\def\S{\mathbb S}
\def\NN{\mathbb N}

\def\CC{\mathbb C}

\def\di{\displaystyle}
\newcommand{\dd}{\,\mathrm{d}}

\def\div{\operatorname{div}}

\def\indicatrice{\textrm{1\kern-0.25emI}} 

\newtheorem{thm}{\bf Theorem}[section]
\newtheorem{lm}[thm]{\bf Lemma}
\newtheorem{prop}[thm]{\bf Proposition}
\newtheorem{cor}[thm]{\bf Corollary}

\newtheorem{rem}[thm]{\bf Remark}

\begin{document}
\title{\bf{Travelling graphs for the forced mean 
curvature motion in an arbitrary space dimension\\
Ondes progressives pour le mouvement par courbure moyenne forc\'e en toute dimension d'espace}}
\author{R\'egis {\sc Monneau}$^{\hbox{\small{a}}}$, Jean-Michel {\sc 
Roquejoffre}$^{\hbox{\small{b}}}$, \\ Violaine  {\sc Roussier-Michon}$^{\hbox{\small{c}}}$\\
\\
\footnotesize{$^{\hbox{a }}$  Universit\'e Paris-Est, Ecole des Ponts ParisTech, CERMICS, 6-8 avenue B. Pascal, Cit\'e
Descartes}\\
\footnotesize{77455 Marne-La-Vall\'ee Cedex 2, France}\\
\footnotesize{$^{\hbox{b }}$ Institut de Math\'ematiques de Toulouse (UMR CNRS 
5219),}\\
\footnotesize{Universit\'e Paul Sabatier, 118 route de Narbonne,  31062 Toulouse 
Cedex 4, France}\\
\footnotesize{$^{\hbox{c }}$ Institut de Math\'ematiques de Toulouse (UMR CNRS 
5219),}\\
\footnotesize{INSA Toulouse, 135 avenue de Rangueil, 31077 Toulouse Cedex 4, 
France}
}
\maketitle

\noindent{\small{{\bf Abstract:} We construct travelling wave graphs of 
the form $z=-ct+\phi(x)$, $\phi: x \in \RR^{N-1} \mapsto \phi(x)\in  \RR$, $N \geq 2$, solutions to the 
$N$-dimensional forced mean curvature motion $V_n=-c_0+\kappa$ ($c\geq c_0$) with prescribed 
asymptotics. For any $1$-homogeneous function $\phi_{\infty}$, viscosity solution to the eikonal 
equation $|D\phi_{\infty}|=\sqrt{(c/c_0)^2-1}$, we exhibit a smooth concave 
solution to the forced mean curvature motion whose asymptotics is driven by 
$\phi_{\infty}$. We also describe $\phi_{\infty}$ in terms of a probability 
measure on $\S^{N-2}$.}} 

\noindent{\small{{\bf R\'esum\'e:} Nous construisons des ondes progressives sous la forme de graphes 
 $z=-ct+\phi(x)$, $\phi: x \in \RR^{N-1} \mapsto \phi(x)\in  \RR$, $N \geq 2$, solutions du mouvement 
par courbure moyenne forc\'e $V_n=-c_0+\kappa$ ($c\geq c_0$) en dimension $N$ d'espace et avec un
comportement asymptotique prescrit. Pour toute solution de viscosit\'e  $\phi_{\infty}$, $1$-homog\`ene en espace,
de l'\'equation eikonale  $|D\phi_{\infty}|=\sqrt{(c/c_0)^2-1}$, nous mettons en \'evidence
une solution r\'eguli\`ere et concave du mouvement par courbure moyenne forc\'e dont le comportement asymptotique 
est donn\'e par $\phi_{\infty}$. Nous d\'ecrivons aussi $\phi_{\infty}$ en terme d'une mesure de probabilité sur la sphère $\S^{N-2}$.}}


\section{Introduction}
\label{intro}

\subsection{Setting of the problem}

The question investigated here is the description of the travelling wave graph
solutions to the forced mean curvature motion in any dimension $N \geq 2$, that is written
under the general form
\begin{equation}
\label{vitesse normale}
V_n= - c_0  + {\kappa}
\end{equation}
where $V_n$ is the normal velocity of the graph, $\kappa$ its local mean curvature 
and $c_0$ a given strictly positive constant to be defined later.
A graph satisfying \eqref{vitesse normale} can be given by the equation 
$z=u(t,x)$ 
where $u:(t,x) \in \RR^+ \times \RR^{N-1} \mapsto u(t,x) \in \RR$ is a solution to 
the 
parabolic equation
\begin{equation}
\label{MCM0}
\frac{u_t}{\sqrt{1+\vert Du\vert^2}}= -c_0 + \mbox{ div }\biggl(\frac{Du}{\sqrt{1+\vert 
Du\vert^2}}\biggl)
 \, , \quad t>0 \, , \, x \in \RR^{N-1}
\end{equation}
Indeed, at any time $t>0$ fixed, the outer normal to the subgraph 
$\{(x,z)\in \RR^{N-1}\times \RR \, | \, z\leq u(t,x)\}$ is given by
$$\vec{n}=\frac{1}{\sqrt{1+|Du|^2}} \left(\begin{array}{c}
-D_x u \\
1
\end{array}\right)$$  
its normal velocity  $V_n$ by $(0,\partial_t u)^T \cdot \vec{n}$ 
while its mean curvature by $\kappa=-\div_{(x,z)} \vec{n}$, see \cite{GT}.

A travelling wave  to \eqref{MCM0} is a solution of the form  $u(t,x)=-ct + 
\phi(x)$ where 
$\phi:x \in \RR^{N-1} \mapsto \phi(x) \in \RR$ is the profile of the wave and 
$c\geq c_0$ is some given constant standing for its speed. 
Thus $\phi$ satisfies the following elliptic equation
\begin{equation}
\label{eq de phi}
- \mbox{ div }\left( \frac{D\phi}{\sqrt{1+|D\phi|^2}}\right) +c_0 - 
\frac{c}{\sqrt{1+|D\phi|^2}}=0 \, , \quad x \in \RR^{N-1}
\end{equation}

\subsection{Connection with reaction diffusion equations}

This work should provide us a better understanding of the multidimensional solutions 
to the  non linear scalar reaction diffusion equation 
\begin{equation}
\label{reaction diffusion}
\partial_t v =\Delta v + f(v) \, , \quad t>0 \, , \, (x,z)\in \RR^{N-1} \times 
\RR
\end{equation}
where $v:(t,x,z) \in [0,+\infty) \times \RR^{N-1} \times \RR \mapsto v(t,x,z) 
\in \RR$ and, especially the case of travelling waves in dimension  $N$.
In the case of a "bistable" nonlinearity $f$, that is to say when $f$ is a 
continuously differentiable function on $\RR$ satisfying 
\begin{description}
\item[\tt i.] $f(0)=f(1)=0$
\item[\tt ii.] $f'(0)<0$ and $f'(1)<0$
\item[\tt iii.] there exists $\theta \in (0,1)$ such that $f(v)<0$ for $v \in 
(0,\theta)$, $f(v)>0$ for $v \in (\theta, 1)$
\item[\tt iv.] $\di\int_0^1 f(v) \dd v >0$,
\end{description}
it is well-known \cite{kanel} that there exists a one-dimensional 
travelling front $v(t,z)=\phi_0(z+c_0 t)$ solution to \eqref{reaction diffusion} 
with $N=1$. 
The speed $c_0$ is unique and strictly positive by $[\tt iv]$ while the profile 
$\phi_0$ is unique up to translations. This result defines the constant $c_0>0$ 
that appears in equation \eqref{vitesse normale}.

In the case $N=2$, multidimensional solutions to \eqref{reaction diffusion} are 
well understood. Paper \cite{HMR1} proves the existence of conical travelling 
waves solutions to \eqref{reaction diffusion}, and paper \cite{HMR3} classifies 
all possible bounded non constant travelling waves solutions under rather weak
 conditions at infinity. In particular, it is  proved in
\cite{HMR3}
that 
$c \geq c_0$ and, up to a shift in $x\in \RR$, either $u$ is a planar front 
$\phi_0(\pm x \cos \alpha+z \sin \alpha)$ with 
$\alpha=\arcsin(c_0/c) \in (0,\frac{\pi}{2}]$
or $u$ is the unique conical front found in \cite{HMR1}. 

In higher dimensions, less is known. In \cite{HMR3}, Hamel, Monneau and 
Roquejoffre proved the existence of conical travelling waves with cylindrical symmetry whose 
level sets are Lipschitz graphs moving away logarithmically from straight cones. 
Some special, non cylindrically symmetric pyramidal-shaped solutions (Taniguchi, \cite{taniguchi2})
are also known in the particular case $N=3$. 

Thus, in order to get a better understanding of the mechanisms at work,
we further the idea of bridging  reaction-diffusion equations with 
geometric motions. In particular, travelling wave graph solutions 
to the forced mean curvature motion go back to Fife \cite{fife}.
He  proved (in a formal fashion) that reaction-diffusion 
travelling fronts propagate with normal velocity
$$V_n=-c_0 + \frac{\kappa}{t} +O\left(\frac{1}{t^2}\right) \, ,  \quad t>>1.$$
For a mathematically rigorous treatment of these ideas, we refer for instance to 
de Mottoni, Schatzman \cite{dMS} - small times, smooth solutions
context - and Barles, Soner, Souganidis \cite{BSS} - arbitrary large times, viscosity solutions context.

Related results must me mentioned in the case of a balanced bistable non-linearity $f$: assumption \texttt{iv.} 
is replaced by $\int_0^1 f(v) \dd v=0$ and \eqref{reaction diffusion} is called the  balanced Allen-Cahn equation. 
In this case, $c_0=0$ and the forced mean curvature equation is replaced by the mean curvature equation. 
Chen, Guo, Hamel, Ninomiya and Roquejoffre \cite{CGHNR} proved that there exist 
cylindrically symmetric traveling waves with paraboloid like interfaces solutions to \eqref{reaction diffusion}  in dimension $N \geq 3$. 
Precisely, they proved that those solutions' level sets are asymptotically given by the equation $z=\frac{c}{2(N-1)} |x|^2$. 
On the other hand, Clutterbuck, Schn\"urer and Schulze \cite{clut} proved that  there exists a unique rotationally symmetric, 
strictly convex, translating graph $u(t,x)=-ct+\phi(r)$ to the mean curvature motion \eqref{eq de phi} 
with $c_0=0$ and whose asymptotics is given by
$$\phi(r)= \frac{c}{2(N-1)} r^2 - \ln r + C + O\left( \frac{1}{r} \right)$$
Further works have also been done in the non radial case for the mean curvature equation. For instance, Xuan Hien Nguyen \cite{nguyen} built non radial and non convex translating graphs solution to \eqref{eq de phi} with $c_0=0$.

\subsection{Main results}

Our theorem \ref{le resultat} below states that, given a $1$-homogeneous solution 
$\phi_\infty$  to the eikonal equation derived from \eqref{eq de phi} (i.e. the equation 
obtained by removing the curvature term)  there exists a smooth solution $\phi$ to  
the forced mean curvature motion equation \eqref{eq de phi} whose asymptotic behaviour is prescribed by
 $\phi_{\infty}$. Here is the precise result.

\begin{thm} 
{\bf (Existence of solutions with prescribed asymptotics in dim. $N$)}
\label{le resultat}\\
Let $N \in \NN \setminus \{0,1\}$, $\alpha \in (0,\frac{\pi}{2}]$, $c_0>0$ and $c=c_0/\sin \alpha$. Choose 
$\phi_{\infty}$ a $1$-homogeneous viscosity solution 
to the  eikonal equation
\begin{equation}
\label{eq de phi infty}
|D\phi_{\infty}(x)|= \cot \alpha\ \, , \quad x\in \RR^{N-1} \, . 
\end{equation}
Then there exists a smooth concave solution $\phi \in C^{\infty}(\RR^{N-1})$ to
\eqref{eq de phi} such that
\begin{equation}
\label{asymptotic of phi}
\phi(x)=\phi_{\infty}(x) + o(|x|) \ \ \hbox{as $\vert x\vert\to+\infty$.}
\end{equation}
\end{thm}
This is the most possible general result. However, due to the possible complexity 
of a solution to the eikonal equation \eqref{eq de phi infty}, it is useful to specialise our result to the particular 
case of a solution with a finite number of facets.

\begin{thm} 
{\bf (Solutions with finite number of facets in dimension $N$)}
\label{le mini resultat}\\
Let $N \in \NN \setminus \{0,1\}$, $\alpha \in (0,\frac{\pi}{2}]$, $c_0>0$ and $c=c_0/\sin \alpha$. 
Choose $\phi^*$ a viscosity solution 
to the  eikonal equation \eqref{eq de phi infty} given for any $x \in \RR^{N-1}$ by
\begin{equation}
\label{eq::profil}
\phi^*(x)=\inf_{\nu \in A} \left( -(\cot\alpha)\  x\cdot \nu + \gamma_\nu \right)
\end{equation}
where $A$ is a finite subset of cardinal $k\in \NN^*$ of the sphere $\S^{N-2}$ and $\gamma_\nu$ are
 given real numbers. Then there exists a  unique smooth concave
 solution $\phi \in C^{\infty}(\RR^{N-1})$ to \eqref{eq de phi} such that
\begin{equation}\label{eq::rr2}
\left\{\begin{array}{l}
\displaystyle -\frac{2 \ln k}{c_0 \sin \alpha}\le \phi-\phi^* \le 0 \, , \quad x \in \RR^{N-1}\\
\\
\displaystyle \lim_{l\to +\infty} \sup_{\mbox{dist}(x,E_\infty)\ge l}|\phi(x)-\phi^*(x)| =0
\end{array}\right.
\end{equation}
where $E_\infty$ is the set of edges defined as
$$E_\infty=\{ x \in \RR^{N-1} \, | \, \phi_\infty  \mbox{ is not $C^1$ at $x$}\}$$
with the  $1$-homogeneous function
$$\phi_\infty(x)=\inf_{\nu \in A} \left( -(\cot\alpha)\  x\cdot \nu \right)$$
\end{thm}

In space dimension $N=3$, we obtain a more precise result by considering solutions 
having a finite number of gradient jumps. Those solutions are still more complex than 
the infimum of a finite number of affine forms.  Here is the precise result. 

\begin{thm} 
{\bf (Solutions with finite number of gradient jumps and $N=3$)}
\label{le micro resultat}\\
Let  $\alpha \in (0,\frac{\pi}{2}]$, $c_0>0$ and $c=c_0/\sin \alpha$. Choose 
$\phi_{\infty}$ a $1$-homogeneous viscosity solution 
to the  eikonal equation \eqref{eq de phi infty} in dimension $N=3$
with a finite number of singularities  on $\S^1$. Then, there exist
\begin{description}
\item[\tt i.] a $2\pi$-periodic continuous function 
$\psi_{\infty}:\theta \in [0,2\pi] \mapsto \psi_{\infty}(\theta) \in [-\cot \alpha, \cot \alpha]$ 
and a finite number $k \in \NN  \backslash \left\{0\right\}$ of  angles $\theta_1<\dots < \theta_k$ in $[0,2\pi)$ such that
$$\phi_{\infty}(r\cos \theta, r \sin \theta)=r \psi_{\infty}(\theta) \, , \quad (r, \theta)  \in \RR^+ \times  [0,2\pi)$$
 Moreover, for any $i \in \{1, \dots, k\}$,
\begin{description}
\item[\tt a.] Either $\forall \theta \in [\theta_i, \theta_{i+1}]$, 
$\psi_{\infty}(\theta)= -(\cot \alpha)$  and we set $\sigma_i=1$
\item[\tt b.] Or
$$\left\{\begin{array}{ll}
\forall \theta \in \left[\theta_i, \frac{\theta_i+ \theta_{i+1}}{2}\right] \, , 
& \psi_{\infty}(\theta)= -(\cot \alpha) \cos(\theta-\theta_i) \\
\forall \theta \in \left[\frac{\theta_i+ \theta_{i+1}}{2}, \theta_{i+1}\right] \, , 
& \psi_{\infty}(\theta)= -(\cot \alpha) \cos(\theta-\theta_{i+1}) 
\end{array}\right.
\mbox{ and we set } \sigma_i=0$$
\end{description}
By convention, $\theta_{k+1}=2\pi + \theta_1$ and $\sigma_{k+1}=\sigma_1$.
If $k\ge 2$, then  $\sigma_i \sigma_{i+1}=0$ for any $i\in \left\{1,...,k\right\}$.
\item[\tt ii.] a smooth concave function $\phi \in C^{\infty}(\RR^2)$ solution  
to equation  \eqref{eq de phi} such that when $|x|$ goes to infinity
$$\phi(x)=\phi_*(x)+O(1)$$
where
\begin{equation}\label{eq::rv22}
\phi_*(x)=-\frac{2}{c_0 \sin \alpha} \ln \left( \int_{\S^1} 
e^{\frac{c_0 \cos \alpha}{2} x \cdot \nu} \dd \mu (\nu) \right)
\end{equation}
and $\mu$ is the non negative measure on $\S^1$ with finite mass  
determined by $\psi_{\infty}$ as follows: 
We set $\mu=\sum_{i=1}^k \mu_i$ where 
for any fixed $\lambda_0>0$, we set
\begin{description}
\item[\tt a.] If $\sigma_i=1$, then $\mu_i= \indicatrice_{(\theta_i,\theta_{i+1})} \dd \theta 
+ \lambda_0 (\delta_{\theta_i} + \delta_{\theta_{i+1}})$\\
(with the exception for $k=1$: $\mu_1=\indicatrice_{(\theta_1,\theta_{1}+2\pi)} \dd \theta$).
\item[\tt b.] If $\sigma_i=0$, then $\mu_i= \lambda_0 (\delta_{\theta_i} + \delta_{\theta_{i+1}})$
\end{description}
\end{description}
\end{thm}

We plan to use  our travelling graphs for the forced mean curvature motion 
exhibited in theorems \ref{le resultat} to \ref{le micro resultat} in order to 
construct  multi-dimensional travelling fronts to the reaction diffusion equation 
\eqref{reaction diffusion}; we plan to do it in a forthcoming paper.

That equation \eqref{eq de phi infty} prescribes the asymptotic 
behaviour of \eqref{eq de phi} has nothing surprising: let 
$\varepsilon >0$ and denote by $\phi_{\varepsilon}$ the scaled function 
$$\phi_{\varepsilon}(x)=\varepsilon \phi\left(\frac{x}{\varepsilon}\right) \, , 
\quad x \in \RR^{N-1}$$
Since $\phi$ is a solution to \eqref{eq de phi}, $\phi_{\varepsilon}$ satisfies 
$$
- \varepsilon \mbox{ div }\left( \frac{D\phi_{\varepsilon}}{\sqrt{1+|D\phi_{\varepsilon}|^2}}\right) 
+c_0 - \frac{c}{\sqrt{1+|D\phi_{\varepsilon}|^2}}=0 \, , \quad x \in \RR^{N-1}
$$
Let $\varepsilon$ go to zero. If adequate estimates for $\phi_{\varepsilon}$ 
are known, (a subsequence of)  $(\phi_\varepsilon)_{\varepsilon>0}$
converges to a function $\phi_{\infty}$ satisfying \eqref{eq de phi infty}.

The proof of Theorem \ref{le resultat} is done by a sub and super solutions
argument.   We first construct a family of
smooth sub-solutions to \eqref{eq de phi}, which will give us some better insight in the 
equation. This step is quite general, and works in any space dimension.
Then, we will construct a Lipschitz super-solution whose 
rescaled asymptotics is prescribed by $\phi_{\infty}$ and this will give us a smooth solution
whose asymptotic behaviour is not well precise. To get a better asymptotics of the
super-solution prescribed by the sub-solution, this will
require a more delicate matching procedure which will limit us, for the moment, 
to any space dimension $N$ with a finite number of facets (theorem \ref{le mini resultat}) or to the
space dimension $N=3$ and a finite number of gradient jumps (theorem \ref{le micro resultat}).

The rest of this paper is organised as follows. In section \ref{section eikonal 
equation}, we build and characterise all $1$-homogeneous solutions to the eikonal 
equation \eqref{eq de phi infty}. In section \ref{section perron}, we detail 
Perron's method in our context, and explain why it will yield a smooth concave solution.
Sub-solutions are built in section \ref{section subsolution}, 
and super-solutions in section \ref{section supersolution}. 
Finally, section \ref{proof le resultat} sums up previous constructions  to prove theorems 
\ref{le resultat} and \ref{le mini resultat}. Section \ref{section sept} presents a more precise 
approach in dimension $N=3$ and details the proof of Theorem \ref{le micro resultat}.
An appendix is devoted to the Laplace's method that we use in our estimates.\\

\noindent
{\bf Acknowledgments.} The first author was partly supported by the ANR project MICA,  
the second and third ones by the  ANR project PREFERED. 
They acknowledge a fruitful discussion with G. Barles and 
 thank C. Imbert for enlightening discussions on his paper \cite{imbert}. They are
also  indebted to H. Berestycki and CAMS Center of EHESS  in Paris for their hospitality 
while preparing this work.

\section{Eikonal equation}
\label{section eikonal equation}
 
In this section, we classify the continuous viscosity solutions to the eikonal 
equation in any dimension $N\geq 2$:
\begin{equation}
\label{eikonal}
|D \phi_{\infty}(x)|=\cot \alpha \, , \quad x \in \RR^{N-1} 
\end{equation}
where $\alpha \in (0,\frac{\pi}{2}]$ is some given angle. 
In a first  subsection, we are interested in the general case. In a second one, we 
reduce our study to $1$-homogeneous functions and give a better description
of those solutions in order to use them in both sections
 \ref{section subsolution} and \ref{section supersolution}.

\subsection{Characterisation of solutions to \eqref{eikonal} in any dimension $N$}
\label{solution eikonal dimension N}

For any unit vector $\nu\in \S^{N-2}$ and $\gamma \in (-\infty,+\infty]$, let us define the affine map
$$\phi_{\nu,\gamma}(x)=-  (\cot \alpha) \ \nu\cdot x + \gamma \in (-\infty,+\infty] \, , \quad x \in \RR^{N-1}$$

\begin{prop}
{\bf (A  Liouville theorem for the eikonal equation)}
\label{th::1}\\
Let $\phi_{\infty} \in C(\RR^{N-1})$. Then $\phi_{\infty}$ is a viscosity solution to the eikonal equation  (\ref{eikonal})
if and only if there exists a lower semi-continuous map $\gamma: \S^{N-2} \to (-\infty,+\infty]$ such that
\begin{equation}
\label{eq::5}
\phi_{\infty}(x)=\inf_{\nu\in \S^{N-2}} \phi_{\nu,\gamma(\nu)}(x)
\end{equation}
Moreover $\phi_{\infty}$ is $1$-homogeneous if and only if   for all $\nu\in\S^{N-2}$, 
$\gamma(\nu)\in \left\{0,+\infty\right\}$.
\end{prop}

This result is most certainly known. Because we not only need the result 
but also an insight of the construction, we give a complete proof.

\noindent {\bf Proof of Proposition \ref{th::1}.}\\
We first show the direct implication. Let $\phi_{\infty} \in C(\RR^{N-1})$ be a viscosity
 solution to \eqref{eikonal}. We shall prove that $\phi_{\infty}$ is $(\cot \alpha)$-Lipschitz 
and concave before giving its characterisation as an infimum of affine maps.

\noindent {\bf Step 1: $\phi_{\infty}$ is locally Lipschitz}\\
Let us consider the ball $B(a,R)$ centered in $a \in \RR^{N-1}$ with radius $R>2$ and define
$$C:=\sup_{|x-y|\le 1,\  (x,y)\in B(a,R)^2} |\phi_{\infty}(x)-\phi_{\infty}(y)|  $$
Because $\phi_{\infty}$ is continuous,  $0 \leq C <+\infty$. 
Denote $\bar{C}=\max(C,\cot\alpha) \in (0,+\infty)$. 
Then we claim that for any $(x,y) \in B(a,R-1)^2$ such that $|x-y| \leq 1$, we have
\begin{equation}\label{eq::1}
|\phi_{\infty}(x)-\phi_{\infty}(y)|\leq \bar{C} |x-y| 
\end{equation}
which asserts that $\phi_{\infty}$ is locally Lipschitz. Indeed, for any point $x_0 \in \overline{B(a,R-1)}$, 
 any constant $\bar{\bar{C}}> \bar{C}$ and  any $\lambda\geq 0$, 
we consider the function $\psi_{\lambda}$ defined as
$$\psi_\lambda(x):=\lambda + \phi_{\infty}(x_0)+ \bar{\bar{C}}|x-x_0| \, , \quad x \in \RR^{N-1}$$
and we set
$$\lambda_*=\inf\left\{\lambda \in \RR^+ \, | \, \forall \mu \geq \lambda \, , \, 
\forall x \in \overline{B(x_0,1)} \, , \quad \psi_\mu(x) \geq \phi_{\infty}(x) \right\}$$
We shall prove by contradiction that $\lambda_*=0$. If not, because $\psi_0 \geq \phi_{\infty}$ 
on $\left\{x_0\right\}\cup \partial B(x_0,1)$, there exists a contact point $z_0$ between 
$\psi_{\lambda_*}$ and $\phi_{\infty}$ which satisfies $z_0\in B(x_0,1)\backslash \left\{x_0\right\}$. 
Then $\psi_{\lambda_*}$ is a test function for the viscosity subsolution $\phi_{\infty}$ at that point. 
Because $|\nabla \psi_{\lambda_*}(z_0)|=  \bar{\bar{C}} > \cot\alpha$, we get a contradiction 
with the viscosity subsolution inequality.
Therefore $\lambda_*=0$ and $\psi_0\geq \phi_{\infty}$ on $B(x_0,1)$. Because this is true for any 
$\bar{\bar{C}}> \bar{C}$, we deduce that this is still true for  
$\bar{\bar{C}}= \bar{C}$ which implies (\ref{eq::1}).\\

\noindent {\bf Step 2: $\phi_{\infty}$ is $(\cot\alpha)$-Lipschitz}\\
We now define 
$$L=\limsup_{n \to +\infty} L_n \quad \mbox{with}\quad 
L_n:= \sup\left\{\frac{\phi_{\infty}(y)-\phi_{\infty}(x)}{|y-x|}, \quad x\in \overline{B(a,R-2)},\quad 
|y-x|\le \frac{1}{n} \right\}$$
Notice that for any $n \in \NN^*$, $0<L_n  \le \bar C$. Moreover, there exists a sequence 
$(x_n,y_n)_{n \in \NN^*}$ such that
$$\lim\limits_{n \to +\infty}\frac{\phi_{\infty}(y_n)-\phi_{\infty}(x_n)}{|y_n-x_n|} = L  
\, \mbox{ and } |y_n-x_n| \leq \frac{1}{n} \,  \mbox{ with } x_n \in \overline{B(a,R-2)}$$
Define for any $x \in \RR^{N-1}$
$$\varepsilon_n=|y_n-x_n| \, , \quad \phi_n(x)=
\frac{\phi_{\infty}(x_n+ \varepsilon_n x)-\phi_{\infty}(x_n)}{\varepsilon_n} 
\quad \mbox{and}\quad \nu_n = \frac{y_n-x_n}{\varepsilon_n} \in \S^{N-2}$$
Thus $\left(\phi_n(\nu_n)\right)_{n \in \NN^*}$ converges to $L$ as $n$ goes to infinity 
and for any $x \in \overline{B(0,1)}$, $|\phi_n(x)| \leq L_n |x|$.
Because $\mbox{Lip}(\phi_n; B(0,n)) \le \bar C$, we see that up to a subsequence,  
$(\phi_n)_{n \in \NN^*}$ converges locally uniformly on $\RR^{N-1}$ to $\phi_0$ 
a viscosity solution to \eqref{eikonal}. Moreover, $\left(\nu_n\right)_{n \in \NN^*}$ 
converges to  $\nu_0\in\S^{N-2}$ with 
$$\phi_0(\nu_0)=L \mbox{ and for any } x \in \RR^{N-1} \, , \quad  \phi_0(x) \le L |x| =:\psi(x)$$
Because $\psi$ touches $\phi_0$ from above at $\nu_0$, we conclude from the 
viscosity inequality for subsolutions that
$$L\le \cot\alpha$$
Now for any $\epsilon>0$, there exists $n_\varepsilon \in \NN^*$ such that 
$L_n \le L +\varepsilon$ for any $n\geq n_\varepsilon$. In particular, for any $(x,y) \in \overline{B(a,R-2)}$
we can split the segment 
$$\displaystyle [x,y] = \bigcup_{i=0,...,K-1} [x_i,x_{i+1}] \mbox{ with } x_0=x  
\, , \, x_{K}=y \mbox{ and } |x_{i+1}-x_i|=\frac{|y-x|}{ K} \le \frac{1}{n}$$ 
This implies that
$$|\phi_{\infty}(x)-\phi_{\infty}(y)|\le (L+\varepsilon) |x-y|$$
which is true for any $\varepsilon>0$. This implies that $\phi$ is $L$-Lipschitz on 
$\overline{B(a,R-2)}$ with $L\le \cot\alpha$.

\noindent {\bf Step 3: $\phi_{\infty}$ is concave}\\
Because  $\phi_{\infty}$  is a Lipschitz stationary viscosity solution to  the evolution equation
$$u_t + H(Du)=0\, ,  \quad x \in \RR^{N-1} \mbox{ where } 
H(p)=\frac12 (p^2- \cot^2\alpha) \, , \quad p \in \RR^{N-1}$$
we can apply Lemma 4 page 131 in \cite{evans}, and get that $\phi_{\infty}$ satisfies for any $t>0$
$$\phi_{\infty}(x+x')-2\phi_{\infty}(x)+\phi_{\infty}(x-x') 
 \leq C_0\frac{|x'|^2}{t},\quad \mbox{for all } (x,x')\in\RR^{2(N-1)}$$ 
and we can check that we have $C_0= 1$.
Letting $t$ go to infinity shows  that $\phi_{\infty}$ is concave in $\RR^{N-1}$.

\noindent {\bf Step 4: Tangent cone}\\
Since $\phi_{\infty}$ is Lipschitz continuous, it is differentiable 
almost everywhere by Rademacher's theorem. 
Let $D \subset \RR^{N-1}$ be the set of differentiability of $\phi_{\infty}$
and fix $x_0\in D$. Since $\phi_{\infty}$ is concave, for any $x \in \RR^{N-1}$, we have 
$$\phi_{\infty}(x) \leq \phi_{\infty}(x_0) + D\phi_{\infty}(x_0) \cdot (x-x_0)$$
Passing to the infimum on $D$, we get  for any $x \in \RR^{N-1}$,
$$\phi_{\infty}(x) \leq \psi(x):= \inf_{x_0 \in D} \phi_{\infty}(x_0) +
 D\phi_{\infty}(x_0) \cdot (x-x_0)$$
Thus, $\psi$ and $\phi_{\infty}$  are $(\cot \alpha)$-Lipschitz functions that coincide on $D$ which
is a dense set on $\RR^{N-1}$. Therefore, they are in fact equal on $\RR^{N-1}$. 
Using equation \eqref{eikonal}, we finally have
$$\phi_{\infty}(x)=\inf_{x_0 \in D} -(\cot \alpha) \, \nu(x_0)\cdot x + g(x_0)$$
where for any $x_0 \in D$, $\nu(x_0)=-D\phi_{\infty}(x_0)/\cot\alpha \in \S^{N-2}$ and 
$g(x_0)=\phi_{\infty}(x_0) - x_0 \cdot D\phi_{\infty}(x_0) \in \RR$. Defining $\gamma$ as
\begin{equation}
\begin{array}{cccl}
\gamma :& \S^{N-2}& \to & (-\infty,+\infty] \\
               & \nu          & \mapsto & \left\{ \begin{array}{l} 
  \inf_{x_0 \in A} g(x_0) \mbox{ if } A:= \{x_0 \in D \, | \, \nu(x_0)=\nu\} \neq \emptyset \\
   +\infty \mbox{ otherwise } \end{array}\right.
\end{array}
\end{equation}
we get the desired characterisation (\ref{eq::5}). 
Since $\phi_{\infty}$ is continuous,
we also deduce from (\ref{eq::5}) that $\gamma$ is lower semi-continuous.

\noindent {\bf Step 5: The $1$-homogeneous case}\\
We assume that $\phi_{\infty}$ is a $1$-homogeneous continuous viscosity solution to \eqref{eikonal}. 
Then for any $x_0\in \RR^{N-1}$,
 there exists $p\in\RR^{N-1}$ with $|p|=\cot\alpha$ such that by \eqref{eq::5}
$$\forall x \in \RR^{N-1} \, , \quad  \phi_{\infty}(x) \leq  \phi_{\infty}(x_0)+p\cdot (x-x_0) $$
On the one hand, considering $x=0$, we get
$$ p\cdot x_0 \leq \phi_{\infty}(x_0)$$
because $\phi_{\infty}$ is $1$-homogeneous. 
On the other hand considering $\lambda x$ instead of $x$ and 
taking the limit $\lambda \to +\infty$, we get
$$\psi(x):=p\cdot x \ge \phi_{\infty}(x) \quad \mbox{with equality at}\quad x=x_0 \, .$$
Therefore if we call ${\mathcal L}_{\phi_{\infty}}$ the set of linear functions $\psi$ 
satisfying $\psi\ge \phi_{\infty}$ such that $|\nabla \psi|=\cot\alpha$,  we have
$$\phi_{\infty} = \inf_{\psi\in {\mathcal L}_{\phi_{\infty}}} \psi$$
because this is true at any point $x_0\in\RR^{N-1}$.\\
\noindent {\bf Step 6: Conclusion}\\
Conversely, if a function $\phi_{\infty}$ is given by (\ref{eq::5}), then it is straightforward 
to check that $\phi_{\infty}$ is a viscosity solution to  (\ref{eikonal}). 
\rule{2mm}{2mm}

\begin{rem} In dimension $N=2$, the previous proposition simply reads:\\
If $N=2$ and $\alpha \in (0,\frac{\pi}{2}]$, $\phi_{\infty}$ is a viscosity 
solution to  \eqref{eikonal} if and only if 
$\phi_{\infty}$ is affine or if there exists $(x_0,y_0) \in \RR^2$ such that
\begin{equation}
\label{valeur absolue}
\phi_{\infty}(x)=- (\cot \alpha) ~ |x-x_0|+y_0 \, , \quad x \in \RR
\end{equation}
Moreover, $\phi_{\infty}$ is $1$-homogeneous if and only if $y_0=0$.

The proof of this proposition can also be done directly from definitions 
of viscosity solutions and we omit the details. 
Notice however the link with \cite{HMR1}: two-dimensional reaction
diffusion waves are either planar fronts or the unique (up to translations) 
conical front whose level sets 
are  asymptotics to the graph of $\phi_{\infty}$ just described.
\end{rem}

\subsection{The $1$-homogeneous case}
\label{homogeneous}

As stressed is theorem \ref{le resultat}, we only build solutions to the forced mean curvature 
motion equation \eqref{eq de phi} whose asymptotics is prescribed by a $1$-homogeneous solution 
to the eikonal equation \eqref{eikonal}. Therefore, it is worth emphasising 
this particular case.

Notice however that there exist viscosity solutions to  
the eikonal equation \eqref{eikonal} defined in $\RR^{N-1}$ that are not homogeneous 
of order $1$.  For instance, consider solutions given by \eqref{valeur absolue}
with $x \in\RR^{N-1}$ and $y_0 \neq 0$. 
We can also consider any translation of a $1$-homogeneous solution.
Another example is for instance given in dimension $N=3$ by a function
$\phi_{\infty}= \inf_{i=1\dots 4} \phi_i$ where $(\phi_i)_{i \in \{1 \dots 4\}}$ 
are four planar solutions defined for $x=(x_1,x_2) \in \RR^2$  by
\begin{align*}
\phi_1(x) = - (\cot \alpha)\  x_1+2 \quad & \quad \phi_2(x)= (\cot \alpha)\  x_1+2 \\
\phi_3(x) = -(\cot\alpha) \ x_2 \quad & \quad \phi_4(x)=(\cot \alpha)\  x_2 
\end{align*}
It is straightforward to check that $\phi_{\infty}$ satisfies 
$|D\phi_{\infty}|=\cot \alpha$ in the viscosity sense and 
that it is not homogeneous of order $1$ since there exists $\lambda >0$ such 
that $\phi_{\infty}(\lambda,0)\neq \lambda \phi_{\infty}(1,0)$.
 
In any case,  a solution $\phi_{\infty}$ to the eikonal equation \eqref{eikonal} 
is concave (see the proof of proposition \ref{th::1},  step 3).  Therefore the function 
$g:  \lambda \in  \RR^{+*} \mapsto g(\lambda)=\phi_{\infty}(\lambda x)/(\lambda |x|) \in \RR$
is decreasing in $\lambda >0$. Since $\phi_{\infty}$ is $(\cot \alpha)$-Lipschitz, 
$g$ is bounded from below and for any 
$x \in \S^{N-2}$, the limit
$$
\lim_{\lambda \to+\infty}\frac{\phi_{\infty}(\lambda x)}{|\lambda|}
$$
exists and $\phi_{\infty}$ is asymptotically homogeneous. 
Thus we have a fairly general understanding
of what is going on by restricting ourselves to homogeneous solutions 
to equation \eqref{eikonal}.

\begin{prop}
\label{countable}
\textbf{(A countable characterisation of homogeneous solutions)}\\
Let $\phi_{\infty} \in C(\RR^{N-1})$. Then $\phi_{\infty}$ is a $1$-homogeneous 
viscosity solution  to the eikonal equation \eqref{eikonal}
if and only if there exists a sequence $(\nu_i)_{i \in \NN}$ of $\S^{N-2}$ such that
\begin{equation}
\label{phi countable}
\phi_{\infty}(x)=\inf_{i \in \NN} -(\cot \alpha) \ \nu_i \cdot x
\end{equation}
\end{prop}

\noindent\textbf{Proof of Proposition \ref{countable}.}\\
Let $\phi_{\infty} \in C(\RR^{N-1})$ be a $1$-homogeneous viscosity solution  to \eqref{eikonal}.
According to proposition \ref{th::1}, there exists a lower semi continuous function
$\gamma$ defined from $\S^{N-2}$ to $\{0,+\infty\}$ such that 
$$\phi_{\infty}(x)=\inf_{\nu \in \S^{N-2}} \phi_{\nu,\gamma(\nu)}(x) \, , \quad x\in \RR^{N-1}$$
Then $K=\{\nu \in \S^{N-2} \, | \, \gamma(\nu) =0\}$ is a compact set of $\S^{N-2}$. 
We claim (see lemma \ref{cube} and corollary \ref{cor cube}  below) 
that there exists  a sequence $(\nu_i)_{i \in \NN}$ of $\S^{N-2}$ such that
$$K=\overline{\bigcup_{i \in \NN} \{\nu_i\}}$$
Thus, $\phi_{\infty}(x)$ can be described as the infimum over $\nu \in K$ of the linear functions
$-(\cot \alpha) \ x \cdot \nu$. Since $\cup_{i \in \NN} \{\nu_i\}$ is dense in $K$, $\phi_{\infty}(x)$
can also be written as the infimum over $i \in \NN$ of the linear functions $-(\cot \alpha) \ x \cdot \nu_i$. 
This ends the proof of proposition \ref{countable} since the converse implication is straightforward. 
\rule{2mm}{2mm}

\begin{lm}
\textbf{(Decomposition of a compact set of $\S^{N-2}$ in cubes)}\\
\label{cube}
For any compact set $K$ of  $\S^{N-2}$, there exists a countable family $(Q_i)_{i \in \NN}$ 
of closed cubes of $\RR^{N-1}$ such that
\begin{equation}
\label{assum}
\left\{\begin{array}{l}
\forall n \in \NN \, , \quad \displaystyle K \subset \bigcup_{i \geq n} Q_i \\
\forall i \in \NN \, , \quad Q_i \cap K \neq \emptyset \\
\displaystyle \limsup_{i \to +\infty} \mbox{diam}(Q_i)=0
\end{array}\right.
\end{equation}
\end{lm}

\noindent\textbf{Proof of Lemma \ref{cube}.}\\
We built this decomposition in cubes by induction. Let $C_0=[-1,1]^{N-1}$ 
be the first cube of width $2$.
Thus  $K \subset C_0$. Since $C_0 \cap K$ is not empty, we divide
$C_0$ in $2^{N-1}$ smaller cubes of width $2^{0}=1$. We call $C_{1,i}$ for 
$i=1 \dots n_1$ those whose intersection 
with $K$ is not empty. Then, $1 \leq n_1 \leq 2^{N-1}$ and 
$$K\quad \subset \quad \bigcup_{i=1}^{n_1} C_{1,i}$$
In the same way, for $i=1,...,n_1$, we divide each cube $C_{1,i}$ in $2^{N-1}$ 
smaller cubes of width $2^{-1}$ and 
keep only those whose intersection with $K$ is not empty. We call them  $C_{2,k}$ 
for $k=1 \dots n_2$ and $1 \leq n_2 \leq  2^{N-1} n_1$. 
Then, one can easily verify that $K \subset \cup_{k=1 \dots n_2} C_{2,k}$.

Assume the cubes $C_{j,i}$ are built for $j\in \NN$, $i=1 \dots n_j$ and $1 \leq n_j \leq 2^{j(N-1)}$
such that 
$$
\left\{\begin{array}{l}
 K \subset \bigcup_{i =1}^{n_j} C_{j,i} \\
\forall i =1 \dots n_j \, , \quad C_{j,i} \cap K \neq \emptyset \\
\mbox{diam}(C_{j,i})  = 2^{-j+1}
\end{array}\right.
$$
Then we construct the cubes $C_{j+1,i}$ as follows. We divide each cube $C_{j,i}$ 
into $2^{N-1}$ smaller cubes of width $2^{-j}$ and 
keep only those whose intersection with $K$ is not empty. We call them $C_{j+1,i}$ 
for $i=1 \dots n_{j+1}$ and $1 \leq n_{j+1} \leq  2^{N-1} n_j\leq 2^{(j+1)(N-1)}$. 
By construction, it is easy to verify that
$$
\left\{\begin{array}{l}
 K \subset \bigcup_{i =1}^{n_{j+1}} C_{j+1,i} \\
\forall i =1 \dots n_{j+1} \, , \quad C_{j+1,i} \cap K \neq \emptyset \\
\mbox{diam}(C_{j+1,i})  = 2^{-j}
\end{array}\right.
$$
The induction is then proved. We thus construct a countable family of cubes that we recall 
$(Q_j)_{j \in \NN}$ for convenience with the desired assumptions \eqref{assum}. 
This ends the proof of lemma \ref{cube}. \rule{2mm}{2mm}

\begin{cor}
\label{cor cube}
\textbf{(Representation of a compact set of $\S^{N-2}$)}\\
For any compact set $K$ of  $\S^{N-2}$, there exists a sequence $(\nu_j)_{j \in \NN}$ 
of $\S^{N-2}$ such that
$$K=\overline{\bigcup_{j \in \NN} \{\nu_j\}}$$
\end{cor}

\noindent \textbf{Proof of Corollary \ref{cor cube}.}\\
For $K$ a compact set of $\S^{N-2}$, we define $(Q_j)_{j \in \NN}$ a family of cubes 
as proposed in lemma \ref{cube}. For any $j \in \NN$, we choose $\nu_j \in K \cap Q_j$.
Then, it is straightforward to check that $\overline{\cup_{j \in\NN} \{\nu_j\}} \subset K$. 
Regarding the converse inclusion, we fix $x_0 \in K$ and $\varepsilon >0$. 
By \eqref{assum}, there exists $n_{\varepsilon} \in \NN$ such 
that the width of cube $Q_i$ is smaller than $\varepsilon$ provided $i \geq n_{\varepsilon}$.
Since $K \subset \cup_{i \geq n_{\varepsilon}}  Q_i$, there exists 
$i_{\varepsilon} \geq n_{\varepsilon}$ such that
$$x_0 \in Q_{i_{\varepsilon}} \, \mbox{ and } \, |x_0-\nu_{i_{\varepsilon}}|\leq \varepsilon  \sqrt{N-1}$$
This shows the density of $\cup_{j \in \NN} \{\nu_j\}$ in $K$ and ends the proof of 
corollary \ref{cor cube}. \rule{2mm}{2mm}

\section{Perron's method and comparison principle}
\label{section perron}

In this section, we are concerned with the forced mean curvature motion equation 
\begin{equation}
\label{MCM}
- \mbox{ div }\left( \frac{D\phi}{\sqrt{1+|D\phi|^2}}\right) +c_0 - 
\frac{c}{\sqrt{1+|D\phi|^2}}=0 \, , \quad x \in \RR^{N-1}
\end{equation}
with the   condition at infinity
\begin{equation}
\label{CL}
\phi(x) =\phi_{\infty}(x) + o(|x|) \, , \quad x \in \RR^{N-1}
\end{equation}
where $\phi_{\infty}$ is a homogeneous viscosity solution to  $|D\phi_{\infty}|=\cot \alpha$
found in section \ref{section eikonal equation} with $\alpha =\arcsin(c_0/c) \in (0,\frac{\pi}{2}]$.
We choose to solve \eqref{MCM} using Perron's method with sub and 
super-solutions (see \cite{user} or \cite{GT}). Let us first recall the existence process and clarify the 
regularity of the solution in the following

\begin{prop}
{\bf (Existence of a  solution to \eqref{MCM} in dimension $N$)}\\
\label{perron}
Let $N \in \NN \setminus \{0,1\}$, $(c_0,c) \in \RR^2$ such that $c \geq c_0>0$.
Assume that $\phi_*$ is a viscosity sub-solution and $\phi^*$  a viscosity super-solution to 
\eqref{MCM} such that $\phi_* \leq \phi^*$ on $\RR^{N-1}$.\\ 
Then,\\
\noindent {\bf i)} there exists a 
function $\phi \in [\phi_*,\phi^*]$ viscosity solution to \eqref{MCM}.\\
\noindent {\bf ii)} Moreover,
if $\phi^*$ is concave, and satisfies the following technical condition:
\begin{equation}
\label{eq::tech}
\mbox{there exists } p\in\RR^{N-1} \mbox{ such that }  \limsup_{|x|\to +\infty} \frac{\phi^*(x)-p\cdot x}{|x|} <0,
\end{equation}
then $\phi$ can be chosen  concave and smooth.
\end{prop}

\noindent\textbf{Proof of Proposition \ref{perron}.}\\ 
We build the solution $\phi$ using Perron's method directly in the framework of 
viscosity solutions to \eqref{MCM}, that is to say $\phi$ is chosen as the maximal 
sub-solution to \eqref{MCM} (see the user's guide to viscosity solutions \cite{user}).\\
\noindent {\bf Step 1: Concavity}\\
We apply a result due to Imbert (see \cite{imbert}) that we first recall. Denote $F$ the following  Hamiltonian 
$$F(p,M)=-\frac{{\mathrm {tr}}M}{\sqrt{1+\vert p\vert^2}}+\frac{{\mathrm {tr}}(M\cdot (p\otimes p))}{(1+\vert p\vert^2)^{3/2}}+c_0-\frac{c}{\sqrt{1+\vert p\vert^2}} \, , \quad (p,M) \in \RR^{N-1} \times \RR^{(N-1)\times (N-1)}_{\textup{sym}}
$$
where $\RR^{(N-1)\times (N-1)}_{\textup{sym}}$ is the set of $(N-1)$-square symmetric matrices.
\begin{prop} 
\label{imbert}
{\bf (Imbert's proposition 5 in \cite{imbert})}\\
 Let $u$ be a lower semi-continuous and epi-pointed function. If $u$ is a supersolution to 
$$F( D u(x), D^2u(x))=0 \, , \quad x \in \RR^{N-1}$$
then, so is its convex envelope.
\end{prop}
In our context, $-\phi^*$ is epi-pointed because of the technical condition (\ref{eq::tech}). 
Thus, the maximal sub-solution $\phi$ to \eqref{MCM} is concave 
(otherwise Imbert's result stated in proposition \ref{imbert} contradicts the 
maximal property of $\phi$). Then, $\phi$ is  a concave viscosity solution to \eqref{MCM}.

\noindent {\bf Step 2: Regularity}\\
Once concavity is at hand, a Lipschitz bound is automatically available from the equation \eqref{MCM}  itself:
$$|D\phi(x)| \leq \cot \alpha \, , \quad x \in \RR^{N-1}$$
where $\alpha \in (0,\frac{\pi}{2}]$ is such that $c_0=c \sin \alpha$.
Then $F$ becomes uniformly elliptic, thus allowing for $C^{1,1}$ estimates (see Theorem 4 in \cite{imbert}).
A bootstrap argument then shows that the solution is $C^\infty$.
This concludes the proof of proposition \ref{perron}.
\rule{2mm}{2mm}

\begin{rem}
Notice that the condition (\ref{eq::tech}) is hidden in the statement of Proposition 9 in \cite{imbert}.
Thus, the proof of  Proposition 5 in \cite{imbert} uses Proposition 9.
\end{rem}
%
%
It now remains to find sub and super-solutions to \eqref{MCM}.

\section{Sub-solution}
\label{section subsolution}

In this section we build smooth sub-solutions to the forced mean curvature equation \eqref{MCM} 
as global solutions to a viscous eikonal equation
and we do believe that they are really close to the desired solutions.

\subsection{Sub-solutions  as solutions to a viscous eikonal equation}

We have the  following
\begin{lm}
\label{CS}
{\bf (Sufficient condition for a sub-solution to \eqref{MCM})}\\
Fix $\alpha \in (0,\frac{\pi}{2}]$, $c_0>0$ and $c=c_0/\sin \alpha$. 
Let $\phi_*$ be a concave smooth solution to 
\begin{equation}
\label{MCM reg concave}
-\Delta \phi_* = \frac{c_0 \sin \alpha}{2} \left(\cot^2 \alpha-|D\phi_*|^2 
\right) \, , \quad x \in \RR^{N-1}
\end{equation}
such that 
\begin{equation}
\label{majoration grad}
|D\phi_*(x)| \leq \cot \alpha \, , \quad x\in \RR^{N-1}
\end{equation}
Then $\phi_*$ is a smooth sub-solution to equation \eqref{MCM}.
\end{lm}

\noindent\textbf{Proof of Lemma \ref{CS}.}\\
Let $\phi_*$ be any concave function verifying 
\eqref{MCM reg concave} and \eqref{majoration grad}. Since $\phi_*$ is 
smooth and concave, we have
\begin{align*}
N[\phi_*]:=- & \mbox{ div }\left( \frac{D\phi_*}{\sqrt{1+|D\phi_*|^2}}\right)  +c_0 - 
\frac{c}{\sqrt{1+|D\phi_*|^2}} \\
                 = & -\frac{\Delta \phi_*}{\sqrt{1+|D\phi_*|^2}} + 
\frac{D^2\phi_*(D\phi_*,D\phi_*)}{(1+|D\phi_*|^2)^{\frac{3}{2}}} +c_0 - 
\frac{c}{\sqrt{1+|D\phi_*|^2}} \\
                 \leq & \frac{1}{\sqrt{1+|D\phi|^2}} \left( -\Delta \phi_* 
+c_0\sqrt{1+|D\phi_*|^2}-c \right)
\end{align*}
From \eqref{majoration grad} and $c_0=c \sin\alpha$, we deduce that 
\begin{align*}
\cot^2\alpha -|D\phi_*|^2& =\left(\frac{c}{c_0}\right)^2 - \left(\sqrt{1+|D\phi_*|^2}\right)^2 \\
         & \leq \frac{2c}{c_0} \left(\frac{c}{c_0}-\sqrt{1+|D\phi_*|^2}\right)
\end{align*}
Using equation \eqref{MCM reg concave} satisfied by $\phi_*$, we get
$$
N[\phi_*] \leq   \frac{1}{\sqrt{1+|D\phi_*|^2}} \left( 
\frac{c}{c_0}-\sqrt{1+|D\phi_*|^2} \right) \left( \frac{2c}{c_0} 
\frac{c_0^2}{2c}-c_0 \right)=0. 
$$
Thus, $\phi_*$ is a sub-solution to \eqref{MCM}.
\rule{2mm}{2mm}\\              

As it is well-known, equation \eqref{MCM reg concave} is readily transformed into
a linear one by the Hopf-Cole transform
$$
\tilde\phi_*(x)={\mathrm{exp}}\biggl(-\frac{c_0\sin\alpha}{2} 
\phi_*\left(\frac{2x}{c_0\cos\alpha}\right)\biggl) \, , \quad x \in \RR^{N-1} 
\, , \quad \alpha \neq \frac{\pi}{2}
$$
where $\tilde \phi_*$ is a positive solution to 
\begin{equation}
\label{e4.1}
-\Delta \tilde{\phi}_*(x) + \tilde{\phi}_*(x) =0 \, , \quad x \in \RR^{N-1}
\end{equation}
From \cite{caf}, a positive solution $\tilde\phi_*$ to \eqref{e4.1} has the form
$$\tilde{\phi}_*(x)= \int_{\S^{N-2}} e^{\nu \cdot x} \dd \mu(\nu) \, , \quad x \in 
\RR^{N-1},
$$
where $\mu$ is a non negative measure on $\S^{N-2}$ with finite mass. 

Now, for any non negative measure $\mu$ on the sphere $\S^{N-2}$, let us define
\begin{equation}
\label{def de psi}
\phi_*(x)= - \frac{2}{c_0 \sin \alpha} \ln \left( \int_{\S^{N-2}} e^{\frac{c_0 
\cos \alpha}{2} x \cdot \nu} \dd \mu(\nu) \right) \, , \quad x \in \RR^{N-1} 
\, , \quad \alpha \in \left(0,\frac{\pi}{2}\right]
\end{equation}
By construction, $\phi_*$ is a smooth solution to \eqref{MCM reg concave}. 
Let us now prove that   $\phi_*$ is a sub-solution to equation \eqref{MCM}, with all the requirements. 

\begin{lm}
\label{lm psi}
{\bf (Inequalities for the derivatives of $\phi_*$)}\\
Let $\mu$  be a non negative measure on $\S^{N-2}$ with finite mass, 
$\alpha \in (0,\frac{\pi}{2}]$, $c_0>0$ and $c=c_0/\sin\alpha$. 
Define $\phi_*$  as in \eqref{def de psi}. 
Then $\phi_*$ is  a smooth concave solution to (\ref{MCM reg concave}) and its gradient is uniformly 
bounded, that is to say, for any $(x,\xi) \in \RR^{N-1} \times \RR^{N-1}$,
\begin{equation}
\label{prop de psi} 
|D\phi_*(x)|\leq \cot \alpha  \, , \quad D^2\phi_*(x)(\xi,\xi) \leq 0 
\end{equation}
\end{lm}

\noindent\textbf{Proof of Lemma \ref{lm psi}}. \\
Let $\mu$ and $\phi_*$ be so defined. We have
$$
D\phi_*(x)=-(\cot \alpha) \, \frac{F_{\nu}(x)}{F_1(x)} \, , \quad x \in \RR^{N-1}
$$
where for any continuous (scalar or vector) function $f$ defined on  $\S^{N-2}$
$$
F_f(x)=\int_{ \S^{N-2}} e^{\frac{c_0 \cos \alpha}{2} x \cdot \nu}f(\nu) \dd\mu(\nu)
$$
Remark that if we define for some fixed $x \in \RR^{N-1}$, 
$\int_{\S^{N-2}} f(\nu) \dd\overline{\mu}_x(\nu)=\frac{F_f(x)}{F_1(x)}$, 
then $\overline{\mu}_x$ is a probability measure on $\S^{N-2}$. 
We can then apply Jensen's inequality to the convex function 
$y \mapsto |y|^2$. This gives
\begin{equation}
\label{Jensen}
\left|\frac{F_{f}(x)}{F_1(x)}\right|^2 =\left| \int_{\S^{N-2}} f(\nu) \dd \overline{\mu}_x \right|^2 
\leq \int_{\S^{N-2}} |f(\nu)|^2 \dd \overline{\mu}_x = \frac{F_{|f|^2}(x)}{F_1(x)}
\end{equation}
for any continuous function $f$ defined on $\S^{N-2}$. Applying this inequality 
to $f(\nu)=\nu$, we get the desired bound on the gradient of $\phi_*$:   
$|D\phi_*(x)|\leq \cot \alpha$. 
 
Regarding the concavity property of $\phi_*$, we use the same type of arguments. 
Indeed, for any $\xi \in \RR^{N-1}$ and $x \in \RR^{N-1}$, we have
$$
D^2\phi_*(x)(\xi,\xi)=-\frac{c_0 \cos^2 \alpha}{2 \sin \alpha} 
\left(\frac{F_{f^2}(x)}{F_1(x)} - \left(\frac{F_f(x)}{F_1(x)}\right)^2\right)
$$
where $f$ is the continuous function defined on $\S^{N-2}$ by $f(\nu)=\nu \cdot \xi$. 
Applying again Jensen's inequality \eqref{Jensen}, we conclude that 
$D^2\phi_*(x)(\xi,\xi) \leq 0$ for any $\xi \in \RR^{N-1}$ and $x \in \RR^{N-1}$ 
which shows that $\phi_*$ is concave.
\rule{2mm}{2mm}

Finally, we proved the following proposition:
\begin{prop}
\textbf{(Existence of a sub-solution to \eqref{MCM})}\\
\label{existence subsolution}
Fix $\alpha \in (0,\frac{\pi}{2}]$, $c_0>0$ and $c=c_0/\sin \alpha$. 
Let $\mu$ be a non negative measure on $\S^{N-2}$ with finite mass.
 Define $\phi_*$ as in \eqref{def de psi}. 
Then, $\phi_*$ is a smooth concave sub-solution to \eqref{MCM}.
\end{prop}

\begin{rem}
The way we choose the measure $\mu$ is decisive in the asymptotic behaviour  
of the sub-solution $\phi_*$
built as in proposition \ref{existence subsolution}. Indeed, if we want the 
subsolution (and hence the solution) 
to the mean curvature equation \eqref{MCM} to follow asymptotically some given 
solution $\phi_{\infty}$ to the eikonal equation \eqref{eikonal}, we will have 
to choose the measure $\mu$ carefully.
In that procedure, information collected in section \ref{section eikonal equation} will help.

Of course, it will be also very interesting to assess whether 
each sub-solution built with a general probability measure gives rise to a solution to 
the mean curvature equation \eqref{MCM}. 
\end{rem}

\section{Super-solution}
\label{section supersolution}

A natural super-solution to the forced mean curvature  equation \eqref{MCM} is a 
viscosity solution $\phi_{\infty}$ to the eikonal equation \eqref{eikonal}. 
Indeed, $\phi_{\infty}$ satisfies (in the distributional and viscosity sense)
$$
-\mbox{ div }\left( \frac{D\phi_{\infty}}{\sqrt{1+|D\phi_{\infty}|^2}}\right) 
+c_0 - \frac{c}{\sqrt{1+|D\phi_{\infty}|^2}} \geq 0 \, , \quad x \in \RR^{N-1}
$$
However, this super-solution does not satisfy the right comparison with the previous 
sub-solution $\phi_*$. For instance, if $N=3$, $\phi_{\infty}$ is the radially symmetric viscosity solution 
to the eikonal equation \eqref{eikonal} and $\phi_*$ the sub-solution  associated with the Lebesgue measure 
$\mu=\dd\theta$ on $[0,2\pi]$ as in \eqref{def de psi}, then we 
can compute the asymptotic behaviour of both functions $\phi_{\infty}$ and 
$\phi_*$ for $x \in \RR^2$ with  $|x|$ large enough 
using Laplace's method (see appendix \ref{laplace}). We then observe that the 
sub-solution $\phi_*$ is above the super-solution $\phi_{\infty}$ in this area. 
This contradicts the crucial assumption $\phi_* \leq \phi^*$ on $\RR^2$ in 
the Perron's method (see proposition \ref{perron}).

\subsection{Super-solutions as infimum of hyperplanes}

Since we do believe that the sub-solution is  close to the viscosity 
solution to the forced mean curvature equation \eqref{MCM} at infinity, 
we prefer to change the super-solution. 
In the general case of dimension $N$, we use the countable 
characterisation of the solution $\phi_{\infty}$  to the eikonal equation 
\eqref{eikonal} that we want to approach (see proposition \ref{countable}). 

\begin{prop}
\label{supersolution dim N}
\textbf{(Existence of a super-solution to \eqref{MCM})}\\
Fix $\alpha \in (0,\frac{\pi}{2}]$, $c_0>0$ and $c=c_0/\sin \alpha$. 
Choose $\phi_{\infty}$ a $1$-homogeneous solution to the eikonal equation \eqref{eikonal}. 
Define $(\nu_i)_{i \in \NN}$ the sequence of $\S^{N-2}$ given by its countable characterisation
in proposition \ref{countable}.

For any sequence $(\lambda_i)_{i \in \NN}$ such that $\lambda_i>0$ and 
$\sum_{i \in \NN} \lambda_i <+\infty$, we set
$$\phi_i(x)=-(\cot \alpha) \ x \cdot \nu_i - \frac{2}{c_0 \sin \alpha} 
\ln \lambda_i \, , \quad i \in \NN \, , \quad x \in \RR^{N-1}$$
 and
$$\phi^*(x)= \inf_{i \in \NN} \phi_i(x) \, , \quad x \in \RR^{N-1}$$
Then, $\phi^*$ is a concave continuous super-solution to \eqref{MCM}.
\end{prop}

\noindent\textbf{Proof of Proposition \ref{supersolution dim N}.}\\
Since $\phi_i$ are exact solutions to the forced mean curvature equation \eqref{MCM},
it is clear that $\phi^*$ is a super-solution to that equation. As the infimum of affine functions, 
it is  concave
and continuous. \rule{2mm}{2mm}

\begin{rem}
This construction is very easy. However, it is not clear whether 
the technical condition \eqref{eq::tech} is satisfied or not. 
It is even clear that when the set of $\{\nu_i\}_{i \in \NN}$ is finite
of cardinal less or equal to $N-1$, this condition is NOT verified. We will see later 
 (see Step 4 of the proof of theorem \ref{le resultat}) how to modify the sub- and 
super-solutions in order to satisfy condition \eqref{eq::tech} and then pass to the limit to recover the
general case.
\end{rem}
\begin{rem}
In the case when the set  $\{\nu_i\}_{i \in \NN}$ is infinite,  the convergence of $\sum \lambda_i$ forces
$(\lambda_i)_{i \in \NN}$ to go to zero and the sequence $(-\ln \lambda_i)_{i \in \NN}$ grows as $i$ goes to infinity.
\end{rem}

\section{ General existence results}
\label{proof le resultat}

Now equipped with sub and super-solutions as well as a Perron's method, 
we are able to prove existence results. 
The general case in dimension $N \geq 2$ is
 the easiest one since the asymptotics is less precise.
Let us explain our ideas in details depending on the degree of precision 
we want to obtain in our construction.

Let $N \in \NN\setminus \{0,1\}$, $\alpha \in (0,\frac{\pi}{2}]$, 
$c_0>0$ and $c=c_0/\sin \alpha$. 
It is worth noticing that some of our constructions do not work 
for $\alpha=\pi/2$. 
However, this case is obvious and leads to planar fronts. Therefore, 
we restrict ourselves to $\alpha \in (0, \frac{\pi}{2})$.

\subsection{Proof of Theorem \ref{le resultat}}
\label{demo main result}

\noindent\textbf{Step 1: Sub and super-solutions}\\
Choose $\phi_{\infty}$ a $1$-homogeneous  continuous viscosity solution 
to the eikonal equation \eqref{eikonal} in $\RR^{N-1}$.
 By proposition \ref{countable}, there exists a 
sequence $(\nu_i)_{i \in \NN}$ of $\S^{N-2}$ such that
$$\phi_{\infty}(x)= \inf_{i \in \NN} -(\cot \alpha) \nu_i \cdot x$$
Let $\mu$ be the probability measure on $\S^{N-2}$  be defined as
$$\mu=\sum_ {i \in \NN} \lambda_i \delta_{\nu_i}$$
where $(\lambda_i)_{i \in \NN}$ are chosen so that $\lambda_i >0$ and 
$\sum_{i =0}^{+\infty} \lambda_i =1$.

Build the sub-solution $\phi_*$ as in \eqref{def de psi} with the 
above measure $\mu$. 
Then by proposition \ref{existence subsolution}, 
$\phi_*$ is a smooth concave sub-solution to the mean curvature 
motion equation \eqref{MCM}. 
Build a concave continuous super-solution $\phi^*$ by 
proposition \ref{supersolution dim N} as the infimum of hyperplanes where the 
$(\lambda_i)_{i \in \NN}$ and $(\nu_i)_{i \in \NN}$ are defined by the choice of $\mu$.
For any $x \in \RR^{N-1}$, 
$$
\phi_*(x) =-\frac{2}{c_0 \sin \alpha} \ln \left( \sum_{i=0}^{+\infty} 
\lambda_i \ e^{\frac{c_0 \cos \alpha}{2} \nu_i \cdot x}\right)
                \leq -\frac{2}{c_0 \sin \alpha} \ln 
\left( \lambda_i \  e^{\frac{c_0 \cos \alpha}{2} \nu_i \cdot x}\right)
$$
Since the last inequality holds for any $i \in \NN$, we have
$$
\phi_*(x) \leq \inf_{i \in \NN } \left(-(\cot \alpha) \nu_i \cdot x -  
\frac{2 }{c_0 \sin \alpha} \ln \lambda_i \right)
= \inf_{i \in \NN}  \phi_i(x)=\phi^*(x)
$$
and the super-solution $\phi^*$  is above the sub-solution $\phi_*$.

\noindent\textbf{Step 2: Asymptotics of sub and super-solutions}\\
Let us now precise their asymptotics: 
we claim that as $|x|$ goes to infinity
\begin{equation}
\label{claim}
\phi_*(x)=\phi_{\infty}(x)+o(|x|) \mbox{ and } \phi^*(x)
=\phi_{\infty}(x) +o(|x|)
\end{equation}
To prove such a claim, the idea is to compare $\phi_{\infty}(x)$ 
with the limits as $\varepsilon$ goes to zero of
$\varepsilon \phi_*(x/\varepsilon)$ and $\varepsilon \phi^*(x/\varepsilon)$. 
In particular, we will prove the sequence of three inequalities: for any $x \in  \RR^{N-1}$ 
\begin{equation}
\label{trois ineq}
\phi_{\infty}(x) \leq \lim\limits_{\varepsilon \to 0} \varepsilon \phi_*(x/\varepsilon) 
\leq \lim\limits_{\varepsilon \to 0} \varepsilon \phi^*(x/\varepsilon) \leq \phi_{\infty}(x)
\end{equation}
which proves the desired claim  \eqref{claim}. 

The first step of the present proof  leads easily to the second inequality in \eqref{trois ineq} 
since $\phi_* \leq \phi^*$ on $\RR^{N-1}$.
As far as the first inequality  is concerned, 
we have for any $i \in \NN$ and $x  \in \RR^{N-1}$,
$$|x \cdot \nu_i| =|x| \cos (\theta_x-\theta_i) \leq |x| \cos \delta_x$$
where $\delta_x$ is the angular distance between $x/|x|$ and 
$K:=\overline{\cup_{i \in \NN} \{\nu_i\}}$ . Thus,
\begin{align*}
\phi_*(x) &\geq -\frac{2}{c_0 \sin \alpha} \ln 
\left( e^{\frac{c_0 \sin \alpha }{2} |x|  \ \cos \delta_x} \mu(\S^{N-2}) \right)\\
               & = - (\cot \alpha) |x|  \ \cos \delta_x= \inf_{\nu \in K} 
-(\cot \alpha) \ x \cdot \nu  =  \phi_{\infty}(x)
\end{align*}
where $\mu(\S^{N-2})=\sum_{i} \lambda_i=1$.
Thus $\phi_{\infty} \leq \phi_*$ on $\RR^{N-1}$ and the homogeneity of 
$\phi_{\infty}$ gives the first inequality of \eqref{trois ineq}.

Regarding the last inequality in \eqref{trois ineq}, 
we know that for any $x \in \RR^{N-1}$,
$$\phi^*(x) \leq \phi_i(x)=-(\cot \alpha) \ x \cdot \nu_i - \frac{2}{c_0 \sin \alpha} \ln \lambda_i$$
Since $\lim_{\varepsilon \to 0} \varepsilon \phi_i(x/\varepsilon)=-(\cot \alpha) \ x \cdot \nu_i$,
 it is clear that
$$\lim_{\varepsilon \to 0}\varepsilon \phi^*\left(\frac{x}{\varepsilon}\right) 
\leq \inf_{i \in\NN} -(\cot \alpha)\ x \cdot \nu_i =\phi_{\infty}(x)$$
This ends the proof of the three inequalities \eqref{trois ineq} and hence of \eqref {claim}.

\noindent\textbf{Step 3: Existence of a solution}\\
By proposition \ref{perron}, there exists a function $\phi \in[\phi_*,\phi^*]$
 viscosity solution to \eqref{MCM}
and by the previous step, $\phi$ verifies the right asymptotics
$$\phi(x)=\phi_{\infty}(x) + o(|x|)$$
However, in the statement of theorem \ref{le resultat}, 
we claim that there exists a smooth concave solution 
to \eqref{MCM} and the above construction does not provide such information. 
By proposition \ref{perron}, the regularity and concavity of the solution
 are at hand if the super-solution $\phi^*$ satisfies 
the technical assumption \eqref{eq::tech}. If it does not,  we will first 
modify the sub and the super-solutions 
in order to satisfy \eqref{eq::tech}, then get a concave solution, and in 
a last step pass to the limit to find a solution (still concave) between $\phi^*$ and $\phi_*$.

\noindent\textbf{Step 4: Regularity and concavity}\\
Let us consider for any $\varepsilon >0$
\begin{equation}
\label{def de psie}
\phi_*^\varepsilon(x)= - \frac{2}{c_0 \sin \alpha} \ln \left( \int_{\S^{N-2}} e^{\frac{c_0 
\cos \alpha}{2} x \cdot \nu} \dd \mu_\varepsilon(\nu) \right) \, , \quad x \in \RR^{N-1}
\end{equation}
with 
$\mu_\varepsilon={\mu} + {\varepsilon} {\mu}_1 $ where
$$ {\mu}_1 = \sum_{\pm}\sum_{j=1}^{N-1} \delta_{\pm e_j}$$
denoting $(e_i)_{i\in \{1,\dots,N-1\}}$ as the canonical 
orthonormal basis of $\RR^{N-1}$.
In the same way, we define 
$$\phi^{\varepsilon*}(x)=\inf_{i \in \NN ,\  j =1 \dots N-1,\ \pm}
\left( -(\cot \alpha)\ x \cdot \nu_i  - \frac{2}{c_0 \sin \alpha} 
\ln {\lambda}_i, -(\cot \alpha)\ x \cdot (\pm e_j) -\frac{2}{c_0 \sin \alpha}\ln {\varepsilon} \right) $$ 
Then, $\phi_*^{\varepsilon}$ is a sub-solution, $\phi^{\varepsilon*}$ is a super-solution
and $\phi_*^{\varepsilon}\le \phi^{\varepsilon*}$. 
It satisfies \eqref{eq::tech} for any $\varepsilon>0$ and for $p=0$.
By proposition \ref{perron}, there exists a concave smooth solution $\phi^\varepsilon$ 
satisfying  equation \eqref{MCM}, with $\phi^\varepsilon$ being $(\cot \alpha)$ - Lipschitz such that
$$\phi_*^{\varepsilon}(x) \leq \phi^{\varepsilon}(x) 
\leq  \phi^{\varepsilon *}(x) \, , \quad x \in \RR^{N-1}$$
Finally, we take the limit as $\varepsilon$ goes to zero. 
The sub-solutions $\phi_*^\varepsilon$ go to $\phi_*$. 
The super-solutions $\phi^{\varepsilon*}$ converge   to $\phi^*$. 
This follows from the expression 
of super-solutions as an infimum of  hyperplanes, 
those associated to the $\varepsilon$ weights going to $+\infty$.
Moreover by Ascoli's theorem, $(\phi^{\varepsilon})_{\varepsilon>0}$ 
converges (up to a subsequence) to some 
concave and $(\cot \alpha)$ - Lipschitz function $\phi^0$ solution to \eqref{MCM} and satisfying
$$\phi_*  \leq \phi^0 \leq  \phi^*$$
Again a bootstrap argument shows that $\phi^0$ is smooth.
Therefore, $\phi^0$ is the intended solution to the mean curvature equation \eqref{MCM}.
\rule{2mm}{2mm}

\subsection{Proof of Theorem \ref{le mini resultat}}

\noindent\textbf{Step 1: Existence of a solution}\\
Choose $\phi^*$ the viscosity solution to the eikonal equation 
\eqref{eikonal}  given by
\begin{equation}
\label{def::31}
\phi^*(x)=\inf_{\nu \in A}( -(\cot\alpha) \ x \cdot \nu + \gamma_{\nu}) \, , \quad x \in \RR^{N-1}
\end{equation}
where $A=\{\nu_1,\dots,\nu_k\}$ is a finite subset of the sphere $\S^{N-2}$, $k \in \NN^*$
and $\gamma_{\nu}$ are any given real numbers.
We build a sub-solution $\phi_*$ as in proposition \ref{existence subsolution}
\begin{equation}
\label{def::32}
\phi_*(x)= -\frac{2}{c_0 \sin\alpha} \ln 
\left( \sum_{i=1}^k \lambda_i \ e^{\frac{c_0 \cos \alpha}{2} x \cdot \nu_i} \right) 
\, , \quad x \in \RR^{N-1}
\end{equation}
where $\lambda_i$ is determined by the relation 
$\gamma_{\nu_i}=-\frac{2}{c_0 \sin \alpha} \ln \lambda_i $ for $i=1 \dots k$.
Let us notice that in the particular case when $A$ is finite, the super-solution
 built in proposition \ref{supersolution dim N} coincides with the solution $\phi^*$
to the eikonal equation.
As in section \ref{demo main result}, $\phi_* \leq \phi^*$ 
and the assumptions of  proposition \ref{perron} i) are satisfied. Thus, there exists 
a function $ \phi \in [\phi_*,\phi^*]$ viscosity solution to \eqref{MCM}. 
Dealing as in section \ref{demo main result} step $4$, we can even find 
a smooth concave solution still denoted $\phi \in [\phi_*,\phi^*]$.
It now remains to study $\phi_* - \phi^*$ to get a precise
asymptotics of the solution $\phi$.\\

\noindent\textbf{Step 2: Asymptotics (first line of (\ref{eq::rr2}))}\\
Setting 
$$ \phi_i(x) =  -(\cot \alpha)\ x \cdot \nu_i  - \frac{2}{c_0 \sin \alpha}  \ln {\lambda}_i$$
we have
\begin{align*}
\phi_*(x) & = -\frac{2}{c_0\sin \alpha}\ln \left(\sum_{i=1}^k e^{-\frac{c_0 \sin \alpha}{2}\phi_i(x)}\right)  
\ge  -\frac{2}{c_0\sin \alpha}\ln\left( k e^{-\frac{c_0 \sin \alpha}{2} 
\left(\displaystyle \min_{i=1,...,k}\phi_i(x)\right)}\right)\\
& = -\frac{2}{c_0\sin \alpha}\ln\left( k e^{-\frac{c_0 \sin \alpha}{2}  \phi^*(x)}\right)
= \phi^*(x) - \frac{2\ln k}{c_0 \sin \alpha}
\end{align*}
This implies in particular that
\begin{equation}\label{eq::rr1}
\displaystyle - \frac{2\ln k}{c_0 \sin \alpha} \le \phi_* - \phi^* \le 0
\end{equation}
which shows the first line of (\ref{eq::rr2}).\\

\noindent\textbf{Step 3: Asymptotics (second line of (\ref{eq::rr2}))}\\
We now notice that the set $E_\infty$ of edges (where $\phi_\infty$ is not $C^1$) is characterized by
$$E_{\infty}=\left\{x\in\RR^{N-1}, \quad \max_{\nu\in A}\  x\cdot \nu = x\cdot \nu_{i_0} = x\cdot \nu_{i_1}, 
 \mbox{ with } \nu_{i_0}\not=\nu_{i_1} \mbox{ and } (\nu_{i_0}, \nu_{i_1}) \in A^2 \right\}$$
For each index $i_0\in \left\{1,...,k\right\}$, let us denote the convex set
$$K_{i_0}=\left\{x\in\RR^{N-1},\quad x\cdot \nu_{i_0} = \max_{\nu\in A}\  x\cdot \nu\right\}$$
Then
$$\partial K_{i_0} \subset \bigcup_{j\not=i_0} (\nu_{i_0}-\nu_j)^{\perp}$$
For $x\in \mbox{Int}(K_{i_0})$, let $x_{i_1}\in \partial K_{i_0}\subset  E_{\infty}$ such that
$$\mbox{dist}(x,E_\infty)= |x-x_{i_1}| \quad \mbox{with}\quad x_{i_1} \in (\nu_{i_0}-\nu_{i_1})^{\perp}.$$
For $j\not= i_0$, we define the orthogonal projection of $x$ on $(\nu_{i_0}-\nu_{j})^{\perp}$ as
$$x_{j}= \mbox{Proj}_{|(\nu_{i_0}-\nu_{j})^{\perp}} (x)$$
In particular $|x-x_j|\ge |x-x_{i_1}|$.
Moreover
\begin{align*}
x\cdot \nu_j & = (x-x_j)\cdot \nu_j + x_j \cdot \nu_j = (x-x_j)\cdot \nu_j + x_j \cdot \nu_{i_0}\\
&= (x-x_j)\cdot (\nu_j-\nu_{i_0}) + x \cdot \nu_{i_0}= x \cdot \nu_{i_0} - |\nu_j-\nu_{i_0}| |x-x_j|\\
& \leq x \cdot \nu_{i_0} - \delta\ \mbox{dist}(x,  E_{\infty})
\end{align*}
with 
$$\delta = \min_{\nu\not=\nu',\ \nu,\nu'\in A} |\nu-\nu'| >0$$
Therefore
\begin{align*}
\phi_*(x)&= \displaystyle  -\frac{2}{c_0\sin\alpha}\ln 
\left(\sum_{i=1}^k\lambda_i e^{\frac{c_0\cos\alpha}{2}\ x\cdot \nu_i}\right)\\
&\ge \displaystyle -\frac{2}{c_0\sin\alpha}\ln 
\left(\sum_{i=1}^k\lambda_i e^{\frac{c_0\cos\alpha}{2}\ 
\left(x\cdot \nu_{i_0} -\delta \ \mbox{dist}(x, E_{\infty})\right)}\right)
\end{align*}
and then for $x\in K_{i_0}$, we have
$$\phi^*(x)\ge \phi_*(x) \ge  \phi^*(x) - 
\frac{2}{c_0\sin\alpha}\ln \left(1+ \sum_{i\not= i_0}\frac{\lambda_i}{\lambda_{i_0}} 
e^{- \frac{c_0\cos\alpha}{2}\ 
\delta\ \mbox{dist}(x,{\color{red} E_\infty})}\right)$$
This shows that
$$\displaystyle \lim_{l\to +\infty} \sup_{\mbox{dist}(x,E_\infty)\ge l}|\phi_*(x)-\phi^*(x)| =0$$
which implies the second line of (\ref{eq::rr2}).

\noindent\textbf{Step 4: Uniqueness}\\
To end the proof of theorem \ref{le mini resultat}, it only remains to prove uniqueness 
of the above smooth solution $\phi$ to the mean curvature equation \eqref{eq de phi} 
with the prescribed asymptotics given by $\phi^*$. Let $\overline{\phi}$ and $\underline{\phi}$ be two 
solutions to \eqref{eq de phi} with the asymptotics (\ref{eq::rr2}). 
Let
$$\varepsilon := \inf \left\{\varepsilon'>0 \, | \, \forall x \in \RR^{N-1} 
\, , \,  \overline{\phi}(x)+\varepsilon' \ge \underline{\phi}(x)\right\}$$
then for any $x \in \RR^{N-1}$, 
$$\overline{\phi}(x)+\varepsilon \ge \underline{\phi}(x)$$
and there exists a sequence of points $(x_n)_n$ such that
$$\overline{\phi}(x_n)+\varepsilon - \underline{\phi}(x_n) \to 0 \,  \mbox{ as $n$ goes to infinity}$$
Let us define for any $x \in \RR^{N-1}$ 
$$\left\{\begin{array}{l}
\overline{\phi}_n(x)=\overline{\phi}(x+x_n)-\overline{\phi}(x_n),\\
\\
\underline{\phi}_n(x)=\underline{\phi}(x+x_n)-\overline{\phi}(x_n)
\end{array}\right.$$
Then, up to the extraction of a subsequence, we have as $n$ goes to infinity
$$
\overline{\phi}_n\to \overline{\phi}_\infty \quad \mbox{ and } \quad 
\underline{\phi}_n \to \underline{\phi}_\infty
$$
with a  uniform convergence on any compact sets of $\RR^{N-1}$.
Moreover $\overline{\phi}_\infty$ and $\underline{\phi}_\infty$ solve equation (\ref{eq de phi})
and satisfy
$$\overline{\phi}_{\infty}+\varepsilon \ge  \underline{\phi}_{\infty} 
\quad \mbox{with equality at}\quad x=0$$
From the strong maximum principle, we deduce that for any $x \in \RR^{N-1}$, 
\begin{equation}\label{eq::rr3}
\overline{\phi}_{\infty}(x)+\varepsilon =  \underline{\phi}_{\infty}(x)
\end{equation}
Let us now assume that $\varepsilon>0$.  Because we have
$$E_\infty\subset \bigcup_{\nu\not= \nu',\ (\nu,\nu')\in A^2} (\nu-\nu')^{\perp} =:\hat{E}_\infty$$
we deduce that there exists $C>0$ such that for any $R\ge 1$ and any $x\in\RR^{N-1}$, we have 
\begin{equation}
\label{eq::dist}
\sup_{y\in \overline{B_R(x)}} \mbox{dist}(y,E_\infty) 
\ge \sup_{y\in \overline{B_R(x)}} \mbox{dist}(y,\hat{E}_\infty) = 
R \sup_{y\in \overline{B_1(x/R)}} \mbox{dist}(y,\hat{E}_\infty)\ge CR
\end{equation}
with
$$C=\inf_{z\in \RR^{N-1}}\left(\sup_{y\in \overline{B_1(z)}} \mbox{dist}(y,\hat{E}_\infty)\right)$$
We easily check by contradiction that $C>0$.
Therefore by \eqref{eq::rr2},  let us choose $R$ large enough such that
$$\displaystyle \sup_{\mbox{dist}(y,E_\infty)\ge CR} |\phi(y)-\phi^*(y)| 
\le \frac{\varepsilon}{4} \quad \mbox{for}\quad \phi=\overline{\phi},\underline{\phi}$$
Then using \eqref{eq::dist},
we get for some $y_n\in\overline{B_R(x_n)}$ with $\mbox{dist}(y_n,E_\infty)\ge CR$,
$$\displaystyle \inf_{y\in\overline{B_R(x_n)}} |\overline{\phi}(y)-\underline{\phi}(y)|
\le |\overline{\phi}(y_n)-\underline{\phi}(y_n)| \le |\overline{\phi}(y_n)-{\phi}^*(y_n)|  
+ |{\phi}^*(y_n)-\underline{\phi}(y_n)| 
\le   \frac{\varepsilon}{2}$$
This implies that
$$\inf_{y\in \overline{B_R(0)}} |\overline{\phi}_\infty(y)-
\underline{\phi}_\infty(y)|\le  \frac{\varepsilon}{2}$$
which is in contradiction with (\ref{eq::rr3}).
Therefore $\varepsilon=0$ and we get $\overline{\phi}\ge \underline{\phi}$.
By symmetry, we also get $\underline{\phi}\ge \overline{\phi}$,
which implies $\overline{\phi} = \underline{\phi}$ and shows the uniqueness of the solution.
This ends the proof of the theorem \ref{le mini resultat}. \rule{2mm}{2mm}

\section{Proof of further results in dimension $N=3$}
\label{section sept}

In this section, the space dimension is $N=3$ and we denote any $x \in \RR^2$ with its polar coordinates
$(r,\theta_x) \in \RR^+ \times [0,2\pi)$ such that $x=r(\cos \theta_x,\sin \theta_x)$.

\subsection{Classification in dimension  $N=3$ of solutions to the eikonal equation with a finite number of singularities}
\label{solution eikonal dimension 3}

This subsection gives  alternative statement and proof of proposition \ref{countable} in dimension $N=3$,
in the special case of a finite number of singularities (i.e. gradient jumps).

\begin{prop}
{\bf (Classification with a finite number of singularities, $N=3$)}
\label{resolution eikonal 3D}\\
Let  $\alpha \in (0,\frac{\pi}{2}]$, $c_0>0$ and $c=c_0/\sin \alpha$. Choose 
$\phi_{\infty}$ a $1$-homogeneous viscosity solution 
to the  eikonal equation \eqref{eikonal} in dimension $N=3$
with a finite number of singularities  on $\S^1$. Then the {\tt i.} of Theorem \ref{le micro resultat} holds.
\end{prop}

\noindent\textbf{Proof of Proposition \ref{resolution eikonal 3D}.}\\
From Proposition \ref{th::1}, we know that there exists a (non empty) compact set $K=\gamma^{-1}(\{0\}) \subset \S^1$, such that
\begin{equation}\label{eq::rv10}
\phi_{\infty}(x)=\inf_{\nu\in K} \left(-(\cot \alpha) \ \nu\cdot x\right)
\end{equation}
Thus, for any $\theta \in [0,2\pi)$, $\psi_{\infty}(\theta)=\phi_\infty(\cos \theta, \sin \theta)$ defines  a continuous function with values in $[-\cot\alpha, \cot \alpha]$. Firstly $\psi_{\infty | K}=-\cot \alpha$. Moreover for any  maximal interval $(a,b)$ contained in $\S^1\backslash K$, we necessarily have
$$\psi_\infty(\theta)=\left\{\begin{array}{ll}
-(\cot \alpha)\ \cos (\theta -a) & \quad \mbox{if}\quad \theta\in \left[a,\frac{a+b}{2}\right],\\
\\
-(\cot \alpha)\ \cos (\theta -b) & \quad \mbox{if}\quad \theta\in \left[\frac{a+b}{2}, b\right].
\end{array}\right.$$
Therefore $\phi_\infty$ has a singularity (gradient jump) at $\theta= \frac{a+b}{2}$. If  $\psi_\infty$ only has a finite number of singularities, then we get the characterization of $\psi_\infty$ given in  the {\tt i.} of Theorem \ref{le micro resultat}.
This ends the proof of proposition \ref{resolution eikonal 3D}.\rule{2mm}{2mm}

\begin{rem}
Notice that without assuming that $\phi_{\infty}$ has a finite number of singularities  on $\S^1$,
the set $K$ could be a Cantor set in (\ref{eq::rv10}).
\end{rem}

\begin{rem}
Notice that the particular function
$\phi_\infty(x) = -(\cot \alpha) |x|$
is the analogue (at the level of the eikonal equation) of the level sets of cylindrically symmetric solutions 
to reaction diffusion equation, constructed in \cite{HMR3} by Hamel, Monneau and Roquejoffre.
Similarly, the particular case where the graph of $\phi_\infty$ 
is a pyramid is also the analogue of solutions contructed by Taniguchi in \cite{taniguchi2}.
\end{rem}

\subsection{Explicit construction of super-solutions in dimension $N=3$}
\label{supersolution dim 3}

In the particular case $N=3$, we construct super-solutions by hand and try to be more precise 
than in section \ref{section supersolution},  above all when $\psi_{\infty}$ is constant and equal 
to $-(\cot \alpha)$ on some interval $I$. In that case, 
we construct our super-solution by hand. We explain  
our ideas on different elementary pieces that we  bring together 
in the proof of theorem \ref{le micro resultat} 
to build a global super-solution $\phi^*$.
Those different elementary pieces are: \textbf{a cone}, \textbf{an edge} or \textbf{an arc} .

\subsubsection{The cone case}

\begin{lm}
\label{supersolution cone}
\textbf{(Radially symmetric solutions)} \\
Let $\phi_{\infty}$ be the viscosity solution to eikonal equation 
\eqref{eikonal} 
whose graph is the straight cone i.e. 
$\phi_{\infty}(x)=-(\cot \alpha) |x|$ for $x \in \RR^2$. 
Then, there exists a unique radially symmetric solution $\phi_c$ (unique up to  an additive constant) 
to the forced mean curvature equation \eqref{MCM},  satisfying 
$$\phi_c'(0)=0 \mbox{ and } \, \phi_c(x)=\phi_{\infty}(x)+o(|x|)$$
Moreover $\phi_c$ is concave and $|D\phi_c|\le \cot \alpha$.
In the case $\alpha=\pi/2$, $\phi_c$ is zero (up to an additive constant). Otherwise,
as $|x|$ goes to infinity, its asymptotics is more precisely given (up to a constant $C \in \RR$)  by
\begin{equation}
\label{asymptotic cone}
\phi_c(x)=-(\cot \alpha) |x| + \frac{1}{c_0 \sin \alpha} \ln |x| 
+ C + \frac{2-3\sin^2 \alpha}{c_0^2 \sin(2\alpha) |x|}+ O \left(\frac{1}{|x|^2}\right) 
\, , \quad \alpha \neq \frac{\pi}{2}
\end{equation}
Moreover, let $\phi_*$ be the sub-solution defined by \eqref{def de psi} with 
$\mu=\frac{d\theta}{2\pi}$ and $N=3$. 
Fix $\phi_c$ such that  
$C=C_0:=\frac{\ln (\pi c_0 \cos \alpha)}{c_0 \sin \alpha}$, 
then for any $x \in \RR^2$, $\phi_c(x) \geq \phi_*(x)$ and as $|x|$ goes to infinity
\begin{equation}\label{eq::rv13}
\phi_c(x)=\phi_*(x)+ O \left(\frac{1}{ \sqrt{|x|}}\right).
\end{equation}
\end{lm}

\noindent\textbf{Proof of Lemma \ref{supersolution cone}.}\\ 
This  result is proved using quite classical methods.
The proof is sketched for the reader's convenience.
With a slight misuse of notation, 
we denote in the case of radially symmetric solutions 
$\phi_c(x)$ by $\phi_c(|x|)=\phi_c(r)$ with $r=|x|\geq 0$. 
Then, equation \eqref{MCM} reads
$$
-\frac{\phi_c'}{r} - \frac{\phi_c"}{1+\phi_c^{'2}} 
+ c_0 \sqrt{1+\phi_c^{'2}} -c=0 \, , \quad r>0
$$
Thus, $\phi_c$ satisfies an ODE involving only its first two derivatives
and  it can only be defined  up to constants. Setting $v=\phi_c'$, we get
\begin{equation}
\label{eq::v}
v'=(1+v^2)\left(c_0 \sqrt{1+v^2}-c-\frac{v}{r}\right):=(1+v^2)g(v,r) \, , \quad r>0
\end{equation}
The proof of lemma \ref{supersolution cone} now reduces to the study of this 
 ODE (existence, uniqueness and asymptotics).

\noindent\textbf{Step 1: Existence} \\
Since for any $r>0$, $g(0,r) \leq 0$ and $g(v_0(r),r)=0$ where 
$$v_0(r)=-\frac{c^2-c_0^2}{\frac{c}{r} + c_0 \sqrt{ \frac{1}{r^2}+c^2-c_0^2}} \leq 0 \, , $$
$v=0$ is a super-solution and $v=v_0$ is a negative decreasing sub-solution to the ODE \eqref{eq::v}.
Thus, for every $r_1>0$, there exist $r_2>r_1$ and 
a solution $v \in C^{\infty}((r_1,r_2),\RR)$ to the ODE $v'=(1+v^2)g(v,r)$ 
satisfying $v_0 \leq v \leq 0$ for any $r \in (r_1,r_2)$. 
Moreover we have $g'_v(v,r) \le 0$ for $v\in [v_0,0]$, and then we conclude that
\begin{equation}\label{eq::rv11}
v'\le 0\quad \mbox{for}\quad r\in(r_1,r_2)
\end{equation}

\noindent\textbf{Step 2: Qualitative properties}  \\
Since for any $r>0$, $v_0(r) \in (-\cot \alpha, 0)$, 
the bounds of $v(r)$ by $v_0(r)$ and zero force $v$ to exist globally for $r>0$. Moreover,
as $\lim\limits_{r \to 0} v_0(r)=0$, $v$ satisfies the same limit  
and we can extend $v$ to $0$ by continuity as $v(0)=0$. 
This proves that $v$ is a global smooth solution to \eqref{eq::v} with initial condition
$v(0)=0$. Thus  it is easy to check that any primitive function $\phi_c$ to $v$ 
is a  smooth radially symmetric solution to \eqref{MCM}
satisfying $\phi_c'(0)=0$.  From (\ref{eq::rv11}), we conclude that $\phi_c$ is concave.

\noindent\textbf{Step 3: Asymptotics} \\
Since $v$ is strictly decreasing on $\RR^+$ and bounded from below, 
it converges to a finite limit $-\cot \alpha \leq l <0$
as $r$ goes to infinity. Since $v$ is uniformly bounded in $[-\cot \alpha,0]$, 
$l$ must satisfy $g(l,+\infty)=0$ which leads to $l=-\cot \alpha$.

Linearising equation \eqref{eq::v} around $-\cot \alpha$,
we set  $w=v+\cot\alpha$. As $w$ is uniformly bounded on $\RR^+$ and 
goes to zero at infinity, equation \eqref{eq::v} reads
$$w'(r)=-c \ (\cot\alpha) w  + \tilde{g}(w,r) \, , \quad  r>0$$
where $\tilde{g}(w,r)=O(w^2) + O(1/r)$ as $r$ goes to infinity. 
By the Duhamel's formula, $w$  follows exponentially fast the 
behaviour of the slowest term of $\tilde{g}$. Thus $w\sim C/r$ as $r$ goes to infinity and 
a straight calculation gives $C=1/(c_0 \sin \alpha)$.
Repeating this method up to order $2$, one gets
$$
v(r)=-\cot \alpha + \frac{1}{c_0 (\sin \alpha) r} 
+ \frac{3\sin^2 \alpha -2}{c_0^2 \sin(2\alpha) r^2}
+ O \left(\frac{1}{r^3}\right)
$$ 
This gives the desired asymptotics for $\phi_c$ up to constants. 

\noindent\textbf{Step 4: Uniqueness}\\
Let $\phi_c^1$ and $\phi_c^2$ be two smooth radially symmetric solutions to \eqref{MCM}. 
From step $3$, we know that they satisfy the same asymptotic expansion
as $r$ goes to infinity and we assume the constants $C$ are the same. 
Since $\phi_c^1-\phi_c^2$ solves an elliptic equation
with smooth co\!efficients and no zero order term, the classical maximum principle applies.
Hence $\phi_c^1-\phi_c^2=0$ because $\lim_{r\to \infty} (\phi_c^1-\phi_c^2)(r)=0$.
This proves the uniqueness of $\phi_c$ up to constants.\\

\noindent\textbf{Step 5: Comparison with $\phi_*$}\\
Using Lemma \ref{l6.1}, we can check (\ref{eq::rv13}) with the suitable value of the constant $C=C_0$
(see for instance the computation (\ref{eq::rv14}) with $N_0(x)\simeq 1/\sqrt{\pi}$).
Finally, using the comparison principle (as in Step 4), we deduce that $\phi_*\le \phi_c$.
This ends the proof of lemma \ref{supersolution cone}.
\rule{2mm}{2mm}

\subsubsection{The edge case}

\begin{lm}
\label{supersolution edge}
\textbf{(Edge super-solution)}\\
Assume $\phi_*$ is given by \eqref{def de psi} 
where the measure $\mu$ is the sum of two Dirac masses
$$
\mu= \mu_{\{\theta_1,\theta_2\}} = \lambda_1 \delta_{\theta_1} + \lambda_2 \delta_{\theta_2}
$$
with $\lambda_i >0$,  $\theta_i \in [0,2\pi)$ for $i=1,2$ such that 
$\theta_1< \theta_2$ and $\delta_{\theta_i}$  the Dirac mass in 
$\theta_i$. 
In the case $\alpha \neq \pi/2$, define $\phi_e$ for any $x \in \RR^2$ by 
\begin{equation}
\label{def phie}
\phi_e(x)=\min(p_1(x), p_2(x)) \mbox{ with } 
p_i(x)= - (\cot \alpha)\  x\cdot \nu_i - \frac{2}{c_0 \sin \alpha} \ln \lambda_i
\end{equation}
where $\nu_i=(\cos \theta_i,\sin \theta_i)$. 
Then, $\phi_e$ is a  Lipschitz and piecewise smooth global super-solution to 
\eqref{MCM} verifying $\phi_* \leq \phi_e$ on $\RR^2$.  
Moreover, as $|x|$ goes to infinity,
\begin{equation}\label{eq::rv12}
\phi_e(x)=\phi_*(x) + O(1) \, , \quad x \in \RR^2
\end{equation}
\end{lm}

\noindent\textbf{Proof of Lemma \ref{supersolution edge}.}\\
Notice first that $\phi_e=\phi^*$ with $\phi^*$ defined as a special case of Proposition \ref{supersolution dim N}.
This shows that $\phi_e$ is a concave (Lipschitz) supersolution. Finally (\ref{eq::rv12}) follows from (\ref{eq::rr1}).
This ends the proof of the lemma.
\rule{2mm}{2mm}

\subsubsection{The arc case}

Here we wish to describe a super-solution to \eqref{MCM} which,
from above, looks like an arc, i.e. is made up of two non parallel straight lines 
 connected by a circle. 

\begin{figure}[h]
\centering\epsfig{figure=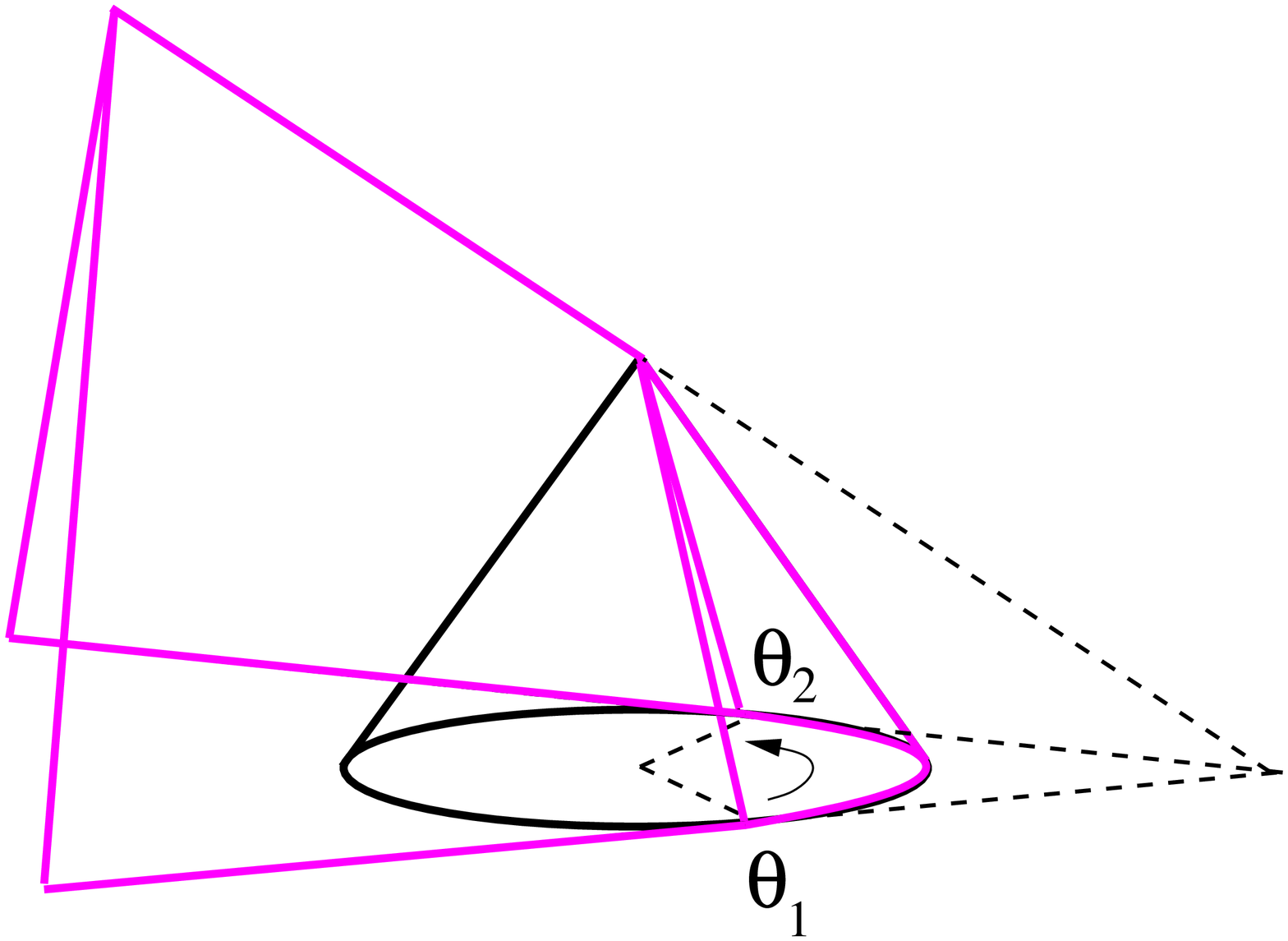,width=80mm}
\caption{Sketch of $\phi_a$ for $\theta_2-\theta_1<\pi$\label{F1}}
\end{figure}

\begin{figure}[h]
\centering\epsfig{figure=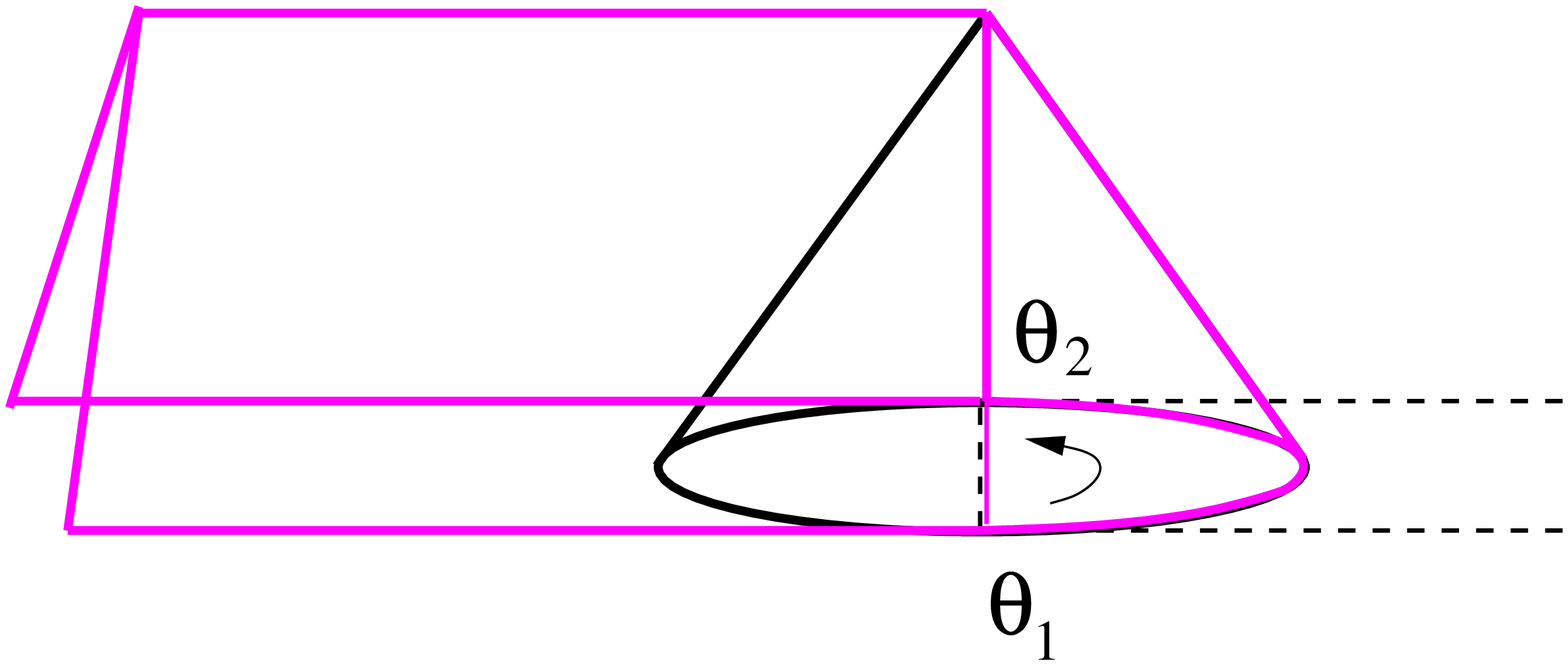,width=80mm}
\caption{Sketch of $\phi_a$ for $\theta_2-\theta_1=\pi$\label{F2}}
\end{figure}

\begin{figure}[h]
\centering\epsfig{figure=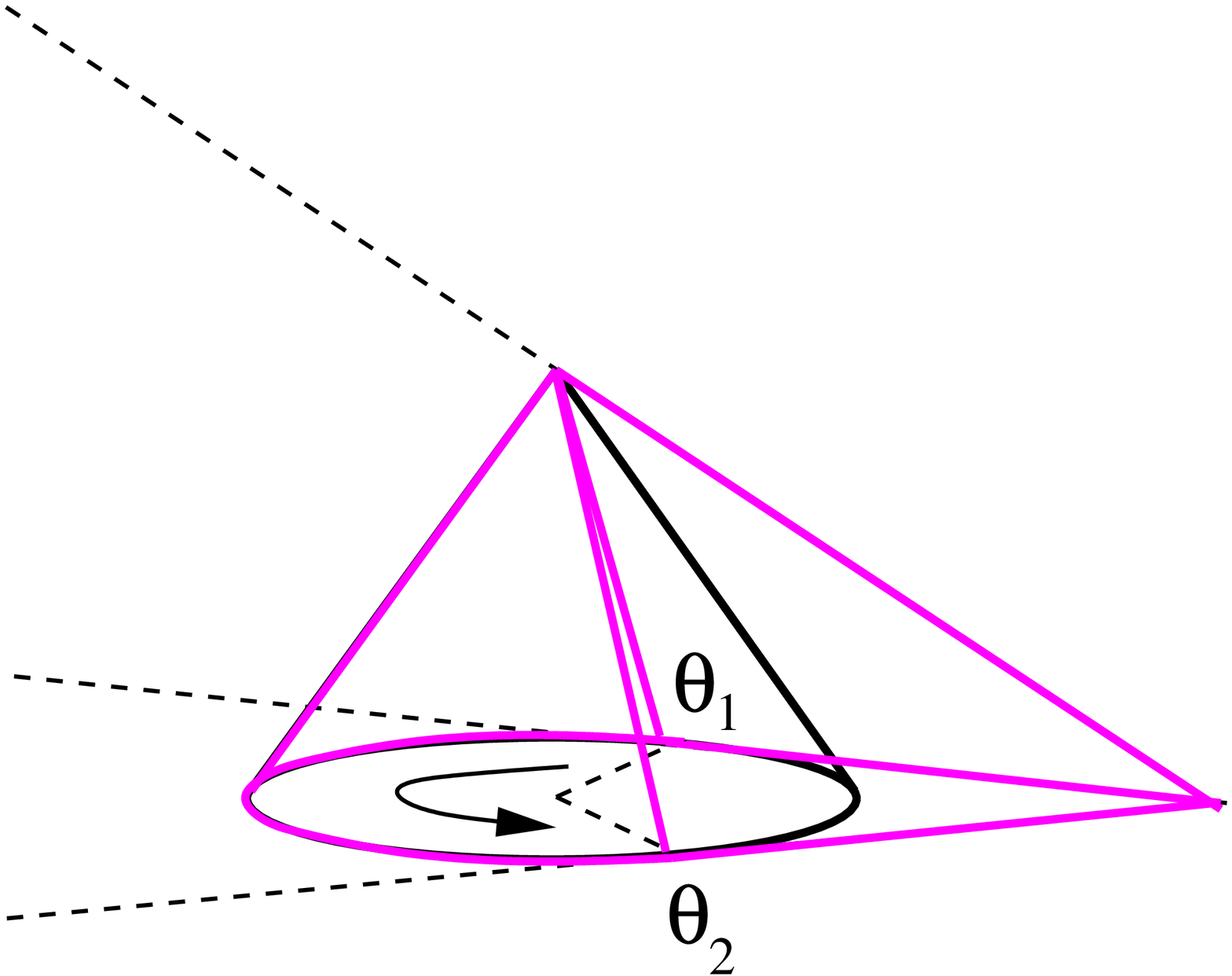,width=80mm}
\caption{Sketch of $\phi_a$ for $\theta_2-\theta_1>\pi$\label{F3}}
\end{figure}

\begin{lm}
\label{supersolution arc} 
\textbf{(Arc super-solution)}\\
Assume $\phi_*$ is given by \eqref{def de psi} 
where the  measure $\mu$ is the sum of two Dirac masses and a Lebesgue measure
$$
\mu= \mu_{[\theta_1,\theta_2]} 
= \lambda \delta_{\theta_1} + \lambda \delta_{\theta_2}
+ \indicatrice_{(\theta_1,\theta_2)}  {\dd \theta},
$$
where $\lambda >0$, $\theta_i \in [0,2\pi)$ for $i=1,2$  and  
$\theta_1<\theta_2$.

Define $\phi_e$ as in the edge case \eqref{def phie} with $\lambda_1=\lambda_2=\lambda$.
Define $\phi_c$ as in the cone case (lemma \ref{supersolution cone})  
where the constant $C \in \RR$ in \eqref{asymptotic cone}  is chosen such that 
$\phi_c(0)=\phi_e(0)= -\frac{2}{c_0 \sin \alpha} \ln \lambda$.

Finally,  define $\phi_a$ on $\RR^2$ by
\begin{equation}
\label{e5.3}
\forall x \in \RR^2 \, , \quad \phi_a(x)=\left\{\begin{array}{ll}
\displaystyle -\frac{2}{c_0 \sin \alpha} \ln \lambda & \mbox{ if } x=0\\
\min (\phi_c(x) , \phi_e(x))& \mbox{  if  }\theta_x \in (\theta_1,\theta_2)\\
\phi_e(x)                   & \mbox{ otherwise, } 
\end{array}\right.
\end{equation}
Then $\phi_a $ is a Lipschitz continuous global super-solution 
to \eqref{MCM}. Moreover,  as $|x|$ goes to infinity 
$$
\phi_a(x)=\phi_*(x)+ O(1)
$$
\end{lm}
The shape of $\phi_a$ is sketched on Figures 
\ref{F1},\ref{F2}, \ref{F3}.\\

\begin{figure}[h]
\centering\epsfig{figure=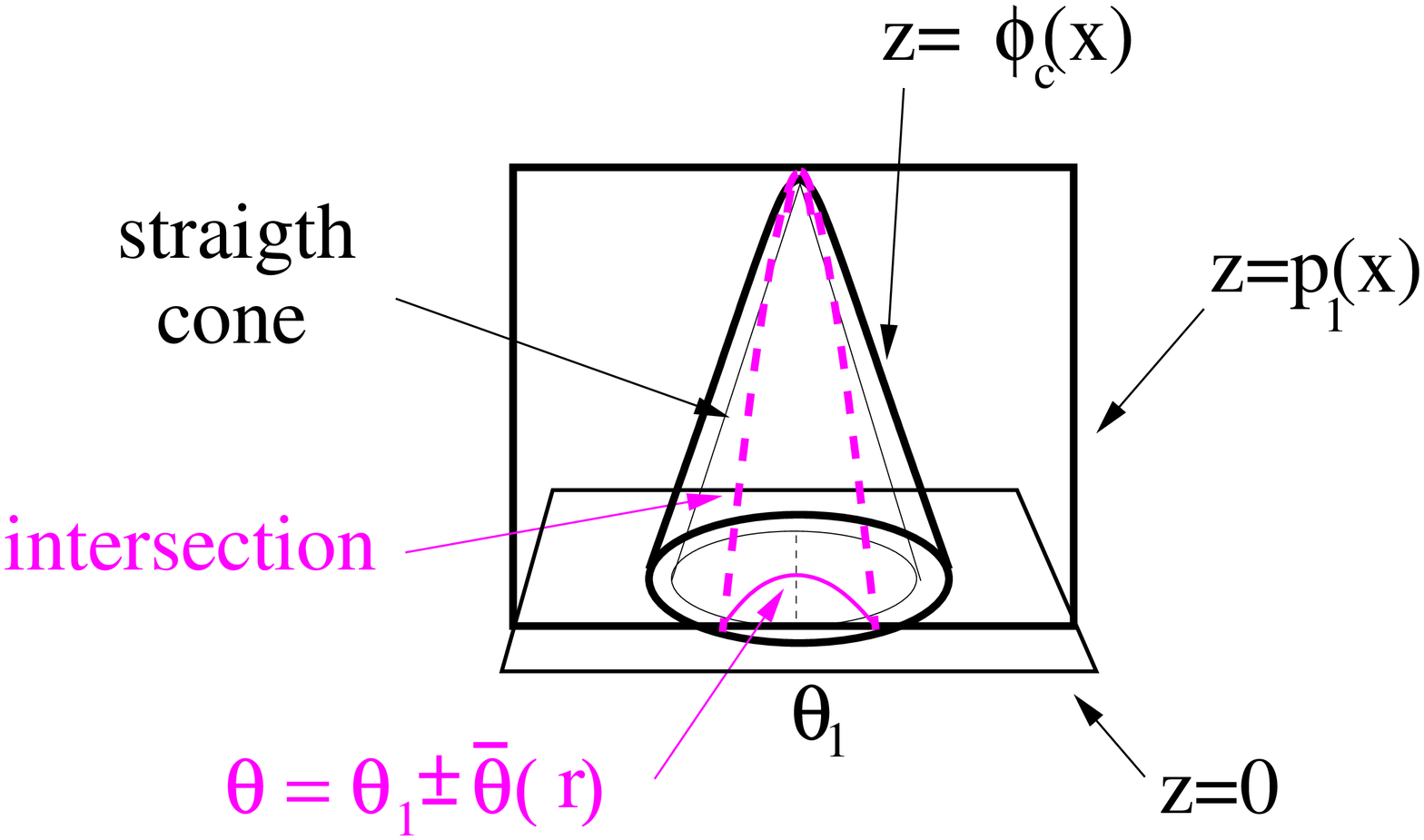,width=120mm}
\caption{Intersection of $z=-(\cot \alpha) \ v_1\cdot x$ with $z=\phi_c(x)$\label{F4}}
\end{figure}

\noindent\textbf{Proof of Lemma \ref{supersolution arc}}\\ 
\noindent\textbf{Step 1: $\phi_a$ is a global continuous super-solution}\\
By definition and lemmas \ref{supersolution cone} and \ref{supersolution edge}, 
$\phi_a$ is a super-solution to \eqref{MCM}
where it is locally the minimum of supersolutions, i.e. everywhere  
except on the two half lines $\theta_x =\theta_i$ for $i=1,2$.
However, we have $\phi_c(0)=\phi_e(0)$ and
$D\phi_e(x)=-(\cot \alpha) \nu_i$ while $\phi_c'(r)\in (- \cot \alpha,0]$ 
for any $x \in \RR^2$ with $\theta_x=\theta_i$, $i=1$ or $2$. Thus,
$\phi_e(x) \leq \phi_c(x)$ on a neighborhood ${\mathcal N}$ (not containing the origin)
of the two half lines $\theta_x = \theta_i$ for $i=1,2$. This implies $\phi_a = \phi_e$ on ${\mathcal N}$ 
and then $\phi_a$ is at least a supersolution on $\RR^2\backslash \left\{0\right\}$.
Moreover, $\phi_a$ is a supersolution on the whole $\RR^2$. Indeed, $\phi_a=\phi_e$ for $\theta_x\not\in [\theta_1,\theta_2]$, then $\phi_a$ has a gradient  jump along the edge $\theta=(\theta_1+\theta_2)/2 +\pi$ up to the origin. 
And this gradient  jump implies that there is no $C^2$ test function touching $\phi_a$ from below at $x=0$.

\begin{figure}[h]
\centering\epsfig{figure=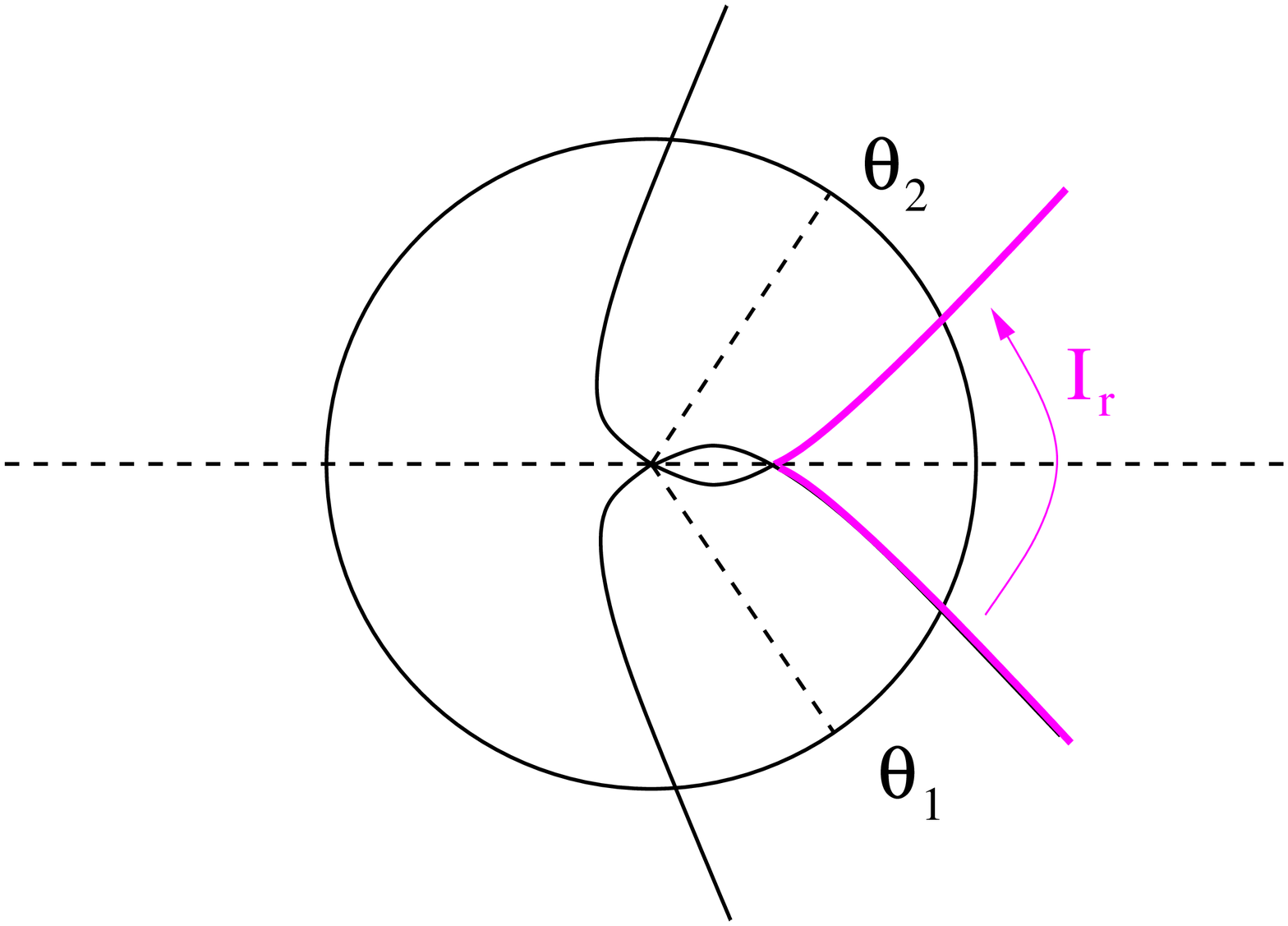,width=80mm}
\caption{The set $I_r$ if $\theta_2-\theta_1<\pi$\label{F5}}
\end{figure}

\begin{figure}[h]
\centering\epsfig{figure=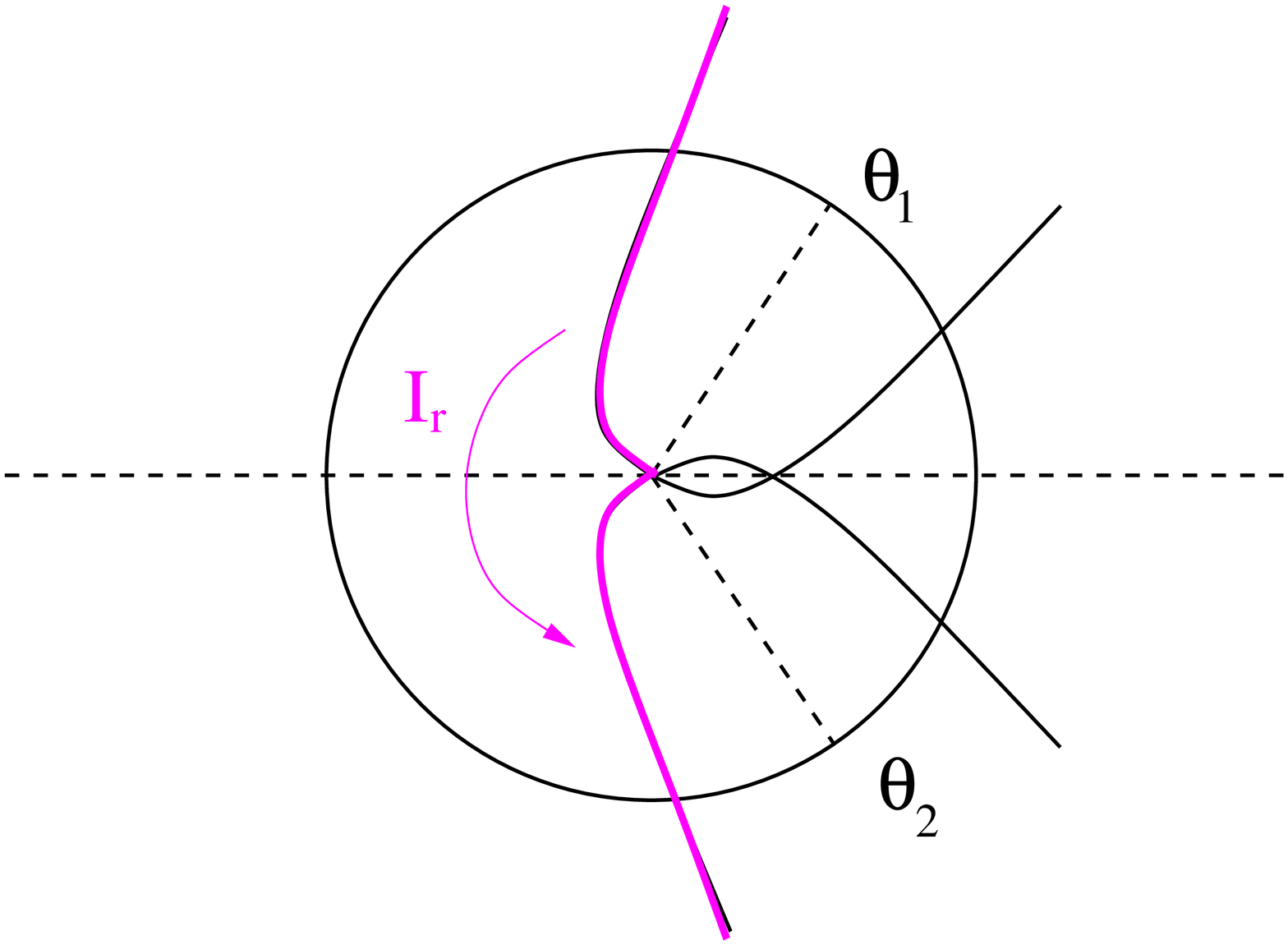,width=80mm}
\caption{The set $I_r$ if $\theta_2-\theta_1>\pi$\label{F6}}
\end{figure}

\noindent\textbf{Step 2: Relative positions of $\phi_e$ and $\phi_c$}\\
Let us now study the relative positions of both graphs of $\phi_e$ and $\phi_c$.
Since $\phi_c(x) \in (-|x| \cot \alpha + \phi_c(0) ,\phi_c(0)]$, we have
\begin{align*}
p_1(x)=\phi_c(x) & \Longleftrightarrow -r (\cot\alpha) \cos( \theta_x-\theta_1)+\phi_c(0)=\phi_c(r) \\
                            & \Longleftrightarrow  \theta_x= \theta_1  \pm 
\arccos\left( \frac{\phi_c(r)-\phi_c(0)}{-r\cot \alpha} \right) := \theta_1 \pm \bar{\theta}(r)
\end{align*}
where $ \bar{\theta}(r) \in (0,\pi/2)$ for $r>0$ (see Figure \ref{F4}). 
Notice that from the concavity of $\phi_c$, we deduce that the set
$$\left\{x,\quad p_1(x) \le \phi_c(x)\right\}=\left\{x,\quad \theta_x\in [\theta_1-\overline{\theta}(r),\theta_1+\overline{\theta}(r)]\right\}$$
is a convex set. Therefore we deduce that
$$ \phi_e(x) \leq \phi_c(x) \Leftrightarrow \theta_x \in 
[\theta_1-\bar{\theta}(r), \theta_1+\bar{\theta}(r)] \cup
  [\theta_2-\bar{\theta}(r), \theta_2+\bar{\theta}(r)]$$
and $\phi_a(x)=\phi_c(x)$ if and only if $\theta_x \in 
I_r= [ \theta_1+\bar{\theta}(r), \theta_2-\bar{\theta}(r)]$.
Since we choose $\phi_c$ such that $\phi_c'(0)=0$, we get
$$\lim_{r \to 0} \bar{\theta}(r)= + \frac{\pi}{2}$$
This forces both curves $\theta_x=\theta_1+\bar{\theta}(r)$ 
and $\theta_x=\theta_2 -\bar{\theta}(r)$ to intersect at some point 
$(x,z) \neq ( 0, \phi_c(0) )$ as soon as $\theta_2-\theta_1 < \pi$. In that case, 
it is worth noticing that the above interval $I_r$ is empty for sufficiently small $r$
(see Figures \ref{F5},\ref{F6}).

On the other hand, using the asymptotics \eqref{asymptotic cone} of $\phi_c$ 
found in lemma \ref{supersolution cone}, we get as $r$ goes to infinity
\begin{equation}
\label{asymptotic cosine}
\cos \bar{\theta}(r) =1- \frac{1}{c_0 \cos \alpha} \frac{\ln r}{r} + 
\frac{\phi_c(0)-C}{r\cot \alpha} + O\left(\frac{1}{r^2}\right) \, , \quad \alpha \neq \frac{\pi}{2}
\end{equation}
where $C$ is the constant given by  \eqref{asymptotic cone} and fixed by the choice $\phi_c(0)=\phi_e(0)$.

\noindent\textbf{Step 3: $\phi_*\le \phi_a + constant$}\\
Let  $\phi_*$ be the sub-solution given by \eqref{def de psi} where 
the  measure $\mu$ is $\mu_{[\theta_1,\theta_2]}$, i.e. with $b =c_0\cos \alpha$:
$$\phi_*(x)=-\frac{2}{c_0\sin \alpha}
\ln\left(
\lambda e^{\frac{br}{2}\cos(\theta_x -\theta_1)} 
+\lambda e^{\frac{br}{2}\cos(\theta_x -\theta_2)} 
+\int_{\theta_1}^{\theta_2}e^{\frac{br}{2}\cos(\theta_x -\theta)}  \ d\theta \right)$$
Each term in the Logarithm being non negative, we have 
\begin{equation}\label{eq::rv16}
\phi_*\leq \phi_e \quad \mbox{on}\quad \RR^2.
\end{equation} 
To prove that $\phi_a$ is above the sub-solution $\phi_*$ up to an additive constant, it remains 
to compare $\phi_* $ and $\phi_c$ when $\theta_x \in 
I_r= [ \theta_1+\bar{\theta}(r), \theta_2-\bar{\theta}(r)]$ 
and $r$ sufficiently large.

According to lemma \ref{l6.1}, one gets that for any $x \in \RR^2$ for $r=|x|$ sufficiently large and 
uniformly in $\theta_x \in [\theta_1,\theta_2]$
\begin{align}\label{eq::rv14}
& \phi_*(x)=-(\cot \alpha) r + \frac{\ln r}{c_0 \sin \alpha} 
- \frac{2}{c_0 \sin \alpha} \ln (\Phi(x))\\
\notag \mbox{ with } \quad  & \Phi(x):=\frac{2\pi N_0(x)}{\sqrt{b}} 
+ \lambda \sqrt{r} e^{\frac{br}{2} (\cos (\theta_x-\theta_1)-1) } 
+ \lambda \sqrt{r} e^{\frac{br}{2} (\cos (\theta_x-\theta_2)-1) } 
+ O\left(\frac{1}{\sqrt{r}}\right)
\end{align}
where 
$$N_0(x)= \int_{\sqrt{r}g(\theta_1-\theta_x)}^{\sqrt{r}g(\theta_2-\theta_x)} 
e^{-\frac{u^2}{4}} \frac{\dd u}{2\pi}$$
as defined in lemma \ref{l6.1}. 
Since $g$ is odd and $\theta_x\in [\theta_1,\theta_2]$,
we see that
$$N_0(x)\ge \int_{0}^{\sqrt{r}g\left(\frac{\theta_2-\theta_1}{2}\right)} 
e^{-\frac{u^2}{4}} \frac{\dd u}{2\pi} = \frac{1}{2\sqrt{\pi}} + o_r(1)$$
We deduce in particular that for $r$ large enough and uniformly in $\theta_x\in[\theta_1,\theta_2]$:
$$\Phi(x)\ge \frac{\sqrt{\pi}}{2\sqrt{b}}$$
Therefore, from the asymptotics (\ref{asymptotic cone}) of $\phi_c$, we deduce that there exist $r_1>0$, $C_1>0$ 
such that  
\begin{equation}
\label{eq::rv15}
\forall r\geq 0 \, , \quad \forall \theta_x \in [\theta_1,\theta_2]\, , \quad r \geq r_1 \Rightarrow \phi_*(x)\le \phi_c(x) + C_1
\end{equation}
Now from (\ref{eq::rv16}) and (\ref{eq::rv15}), we deduce that (up to increasing the constant $C_1$),
$$\forall r\geq 0 \, , \quad \forall \theta_x \in [\theta_1,\theta_2]\, , \quad r \geq r_1 \Rightarrow \phi_*(x)\le \phi_a(x) + C_1$$

\noindent\textbf{Step 4: $\phi_*\ge \phi_a - constant$}\\

\noindent\textbf{Case 1: $\theta_x\in I_r$}\\
We start with the asymptotics (\ref{eq::rv14}).
Using (\ref{asymptotic cosine}), we see that there exist $r_2>0$, $C_2>0$ such that 
$$\forall i=1,2 \, , \quad \forall r\geq 0 \, , \quad \forall \theta_x \in I_r \, , \quad 
r \geq r_2 \Rightarrow\sqrt{r} e^{\frac{br}{2} (\cos (\theta_x-\theta_i)-1)} 
\le \sqrt{r} e^{\frac{br}{2} (\cos (\overline{\theta}(r))-1)}\le C_2$$
Using also the fact that $N_0(x)\le 1/(2\sqrt{\pi})$, we deduce that
$\Phi$ is bounded for $r$ large enough and then (up to increasing $r_2$ and $C_2$)
\begin{equation}
\label{eq::rv17}
 \forall r\geq 0 \, , \quad \forall \theta_x \in I_r \, , \quad 
r \geq r_2 \Rightarrow \phi_*(x) \ge \phi_c(x) -C_2
\end{equation}

\noindent\textbf{Case 2: $\theta_x\in [\theta_1,\theta_2]\backslash I_r$}\\
Let us assume that $\theta_x \in [\theta_1,\theta_1+\overline{\theta}(r))$ (the symmetric case is similar).
Then there exist $r_3>0$, $C_3>0$ such that
$$\forall r \geq 0 \, , \quad \forall \theta_x \in [\theta_1,\theta_1+\overline{\theta}(r)) \, , \quad r \geq r_3 \Rightarrow
\sqrt{r} e^{\frac{br}{2} (\cos (\theta_x-\theta_1)-1)}\ge  \sqrt{r} e^{\frac{br}{2} (\cos (\overline{\theta}(r))-1)} \ge C_3>0$$
Therefore,( up to increasing the constants $r_3$ and $C_3$) for any $r\geq 0$ and any $\theta_x \in [\theta_1,\theta_1+\overline{\theta}(r))$,
$$r \geq r_3 \Rightarrow \Phi(x)\le C_3 \left(\lambda \sqrt{r} e^{\frac{br}{2} (\cos (\theta_x-\theta_1)-1) } 
+ \lambda \sqrt{r} e^{\frac{br}{2} (\cos (\theta_x-\theta_2)-1) }\right)$$
and then
\begin{align}
\notag 
\phi_*(x) & \ge \displaystyle - \frac{2}{c_0 \sin \alpha} \ln \left\{ \lambda e^{\frac{br}{2}\cos(\theta_x -\theta_1)} 
+\lambda e^{\frac{br}{2}\cos(\theta_x -\theta_2)} \right\} - \frac{2\ln C_3}{c_0 \sin \alpha}\\
\notag & \geq -\frac{2}{c_0 \sin \alpha} \ln \left( 2 e^{\frac{c_0 \sin \alpha}{2} \phi_e(x) }\right)  - \frac{2\ln C_3}{c_0 \sin \alpha}\\
\label{eq::rv20} & = \displaystyle \phi_e(x) - \frac{2\ln (2C_3)}{c_0 \sin \alpha}
\end{align}

\noindent\textbf{Case 3: $\theta_x\in \S^1\backslash [\theta_1,\theta_2]$}\\
Notice that the set $\S^1\backslash [\theta_1,\theta_2]$ is not empty because $\theta_2-\theta_1< 2\pi$ (as a consequence of $\theta_1,\theta_2\in [0,2\pi)$).
In that case, we define 
$$\theta_m'=\theta_m -\pi\quad \mbox{with}\quad \theta_m=\frac{\theta_1+\theta_2}{2}.$$
Then, it satisfies 
$\theta_2-2\pi < \theta_m' < \theta_1$. Let us assume that (the other case is similar):
\begin{equation}\label{eq::rv19}
\theta_x\in [\theta_m',\theta_1).
\end{equation}
We also define
$$\theta_x'= \theta_x+\pi$$
Then we have
$$\int_{\theta_1}^{\theta_2}e^{\frac{br}{2}\cos(\theta_x -\theta)}  \ d\theta  
 = \int_{[\theta_1,\theta_2]\cap [\theta_x,\theta_x']} (...)\ d\theta + \int_{[\theta_1,\theta_2]\backslash [\theta_x,\theta_x']} (...)\ d\theta$$
We have
$$\int_{[\theta_1,\theta_2]\cap [\theta_x,\theta_x']}e^{\frac{br}{2}\cos(\theta_x -\theta)}  \ d\theta \le \pi \ e^{\frac{br}{2}\cos(\theta_x -\theta_1)}$$
Using (\ref{eq::rv19}), we also see that
$$\int_{[\theta_1,\theta_2]\backslash  [\theta_x,\theta_x']}e^{\frac{br}{2}\cos(\theta_x -\theta)}  \ d\theta \le \pi \ e^{\frac{br}{2}\cos(\theta_x -\theta_1)}$$
Therefore, we conclude that in this third case,
\begin{align}
\notag \phi_*(x)
& \ge \displaystyle - \frac{2}{c_0 \sin \alpha} \ln \left\{ (\lambda+2\pi) e^{\frac{br}{2}\cos(\theta_x -\theta_1)} 
+(\lambda+2\pi) e^{\frac{br}{2}\cos(\theta_x -\theta_2)} \right\}\\
\label{eq::rv21} & \ge  \displaystyle \phi_e(x) - \frac{2\ln \left(\frac{2(\lambda+2\pi)}{\lambda}\right)}{c_0 \sin \alpha} 
\end{align}

\noindent\textbf{Conclusion}\\
Putting (\ref{eq::rv17}), (\ref{eq::rv20}) and (\ref{eq::rv21}) together, we get that there exists a constant $C>0$ such that 
for $r$ large enough and uniformly in $\theta_x\in \S^1$.
$$\phi_*(x) \ge \phi_a(x) -C$$
The functions $\phi_*$ and $\phi_a$ being continuous, the result still holds for any $r \geq 0$ (up to increasing the constant $C$).
This concludes the proof of lemma \ref{supersolution arc}.
\rule{2mm}{2mm}

\subsection{Proof of Theorem \ref{le micro resultat}}

\noindent\textbf{Step 1: Existence of  a solution}\\
Let  $\alpha \in (0,\frac{\pi}{2}]$, $c_0>0$ and $c=c_0/\sin \alpha$. 
The case $\alpha=\pi/2$ is obvious and we omit it.
Choose $\phi_{\infty}$ a $1$-homogeneous viscosity solution 
to the  eikonal equation (\ref{eikonal}) in dimension $N=3$ with a finite number $m$ of singularities.
By proposition \ref{resolution eikonal 3D}, the {\tt i.} of Theorem \ref{le micro resultat} is already established,
and we can consider the measure $\mu$ given in the {\tt ii.} of Theorem \ref{le micro resultat}.
Then Proposition \ref{existence subsolution} 
implies that the function $\phi_*$ given by (\ref{eq::rv22}) 
is a smooth concave subsolution of (\ref{MCM}).

If $k=1$, $\phi_{\infty}$ has no gradient jump and the corresponding measure is $\mu=\dd \theta$ 
or $\mu=\lambda_0 (\delta_{\theta_1}+\delta_{\theta_1+2\pi})$. In the first case, we saw in lemma \ref{supersolution cone} that $\phi_c$, 
to which a suitable constant is added, is a smooth solution to \eqref{MCM} with the right asymptotics at infinity.
In the second one, $\phi=\phi_{\infty}$ is a suitable solution to \eqref{MCM}.

We now turn to the case $k \ge 2$. 
For any $i \in \{1, \dots ,k\}$, choose $\lambda_0>0$ a given positive constant.
We have
$$\mu=\sum_j \mu_j \ge \tilde{\mu}_i$$
with
\begin{equation}
\label{def de mu i tilde}
\tilde{\mu}_i=
2\lambda_0 (\delta_{\theta_i} + \delta_{\theta_{i+1}}) +\sigma_i\ \indicatrice_{(\theta_i,\theta_{i+1})} \dd \theta
\end{equation}
Let $\tilde{\phi}_{i*}$ be the subsolution defined in \eqref{def de psi} with the measure  $\tilde{\mu}_i$.
If $\tilde{\mu}_i$ corresponds to an arc ($\sigma_i=1$), denote $\tilde{\phi}_i^*$ the global supersolution 
defined in lemma \ref{supersolution arc} with $\lambda=2\lambda_0$. Notice that there is a constant $C>0$ (that can be chosen independently of the index $i$) such that
\begin{equation}\label{eq::rv25}
\tilde{\phi}^*_{i} -C\le \tilde{\phi}_{i*} \le \tilde{\phi}^*_{i} +C
\end{equation}
If $\tilde{\mu}_i$ corresponds to an edge ($\sigma_i=0$), denote $\tilde{\phi}_i^*$ the global supersolution 
defined in lemma \ref{supersolution edge} with $\lambda_1=\lambda_2=2\lambda_0$, 
which satisfies in particular (\ref{eq::rv25}).
Finally, define on $\RR^2$ the function $\tilde{\phi}^*$ as the infimum over $i \in \{1, \dots ,k\} $ of $\tilde{\phi}_i^*$.
Notice that, by construction, we have
\begin{equation}\label{eq::rv40}
\tilde{\phi}^*(x)=\tilde{\phi}_i^*(x) \quad \mbox{if}\quad \theta_x\in [\theta_i,\theta_{i+1}]
\end{equation}
We also have in particular
$$\phi_*\le \tilde{\phi}_{i*} \le \tilde{\phi}^*_{i} +C \le \tilde{\phi}^* +C =: \phi^*$$
We claim that at infinity 
\begin{equation}\label{eq::rv26}
\phi^*(x) = \phi_*(x) +O(1). 
\end{equation}
We shall first finish the proof of theorem \ref{le micro resultat} 
and come back to the proof of that claim in a second step.
Thus $\phi^*$ is a global supersolution above the subsolution $\phi_*$.
Moreover, 
either there exists $\sigma_i=1$ and then we have (see in particular Figures \ref{F1}, \ref{F2}, \ref{F3}) 
\begin{equation}
\label{eq::rv45}
\mbox{there exists } p\in\RR^{N-1} \mbox{ such that }  \limsup_{|x|\to +\infty} \frac{\phi^*(x)-p\cdot x}{|x|} <0.
\end{equation}
Or $\sigma_i=0$ for any $i$, and condition (\ref{eq::rv45}) is satisfied if $k\ge 3$.
The special case $k=2$ and $\sigma_1=\sigma_2=0$ corresponds to an edge for which we already know the existence of a smooth concave solution, by Theorem \ref{le mini resultat}.
In the other cases, condition (\ref{eq::rv45}) and Proposition \ref{perron}
imply the existence of a smooth concave solution $\phi \in [\phi_*,\phi^*]$.

\noindent\textbf{Step 2: Proof of (\ref{eq::rv26}) in the case $k \ge 2$}\\
Let $x \in \RR^2$, then there exists $i \in \{1, \dots, k \}$ such that 
$\theta_x \in [\theta_i,\theta_{i+1}]$.\\
We can write:
$$\mu = \tilde{\mu}_i +\bar\mu_i$$
where $\tilde{\mu}_i$ is defined by \eqref{def de mu i tilde}. So
\begin{equation}
\label{eq::rv31}
\mbox{supp}(\bar\mu_i)\quad  \subset \quad \S^1\backslash (\theta_i,\theta_{i+1}) \quad \not=\emptyset
\end{equation}
Then $k\ge 2$ implies that
$$l_x:=\min\left\{\theta_x -\theta_i, \theta_{i+1}-\theta_x\right\}< \pi$$
and
$$\int_{\S^1}e^{\frac{br}{2}\cos(\theta_x-\theta)}\ d\bar \mu_i(\theta) 
\le \bar \mu_i(\S^1) e^{\frac{br}{2}\cos l_x}\le  \frac{\bar \mu_i(\S^1)}{2\lambda_0}\
\int_{\S^1}e^{\frac{br}{2}\cos(\theta_x-\theta)}\ d\tilde{\mu}_i(\theta)$$
with
$$\int_{\S^1}e^{\frac{br}{2}\cos(\theta_x-\theta)}\ d\tilde{\mu}_i(\theta) = 
2\lambda_0 e^{ \frac{br}{2} \cos (\theta_x-\theta_i)}  
 + 2\lambda_0 e^{ \frac{br}{2} \cos (\theta_x-\theta_{i+1})} 
 + \sigma_i \int_{\theta_i}^{\theta_{i+1}} e^{ \frac{br}{2} \cos (\theta_x-\theta)} \dd \theta$$
Therefore we have
\begin{align*}
\phi_*(x)& \ge -\frac{2}{c_0 \sin \alpha} \ln \left\{\left(1+\frac{\bar \mu_i(\S^1)}{2\lambda_0}\right)
\int_{\S^1}e^{\frac{br}{2}\cos(\theta_x-\theta)}\ d\tilde{\mu}_i(\theta)\right\}\\
& =  \tilde{\phi}_{i*}(x)  -\frac{2}{c_0 \sin \alpha} \ln \left(1+\frac{\bar \mu_i(\S^1)}{2\lambda_0}\right)\\
& \ge  \tilde{\phi}^*_{i}(x) -C'
\end{align*}
where we have used (\ref{eq::rv25}) in the last line.
Using (\ref{eq::rv40}), we see that this implies
\begin{equation}\label{eq::rv41}
\phi_*(x) \ge {\phi}^*(x) -C'' \quad \mbox{for}\quad \theta_x\in [\theta_i,\theta_{i+1}] 
\end{equation}
Finally, this implies (\ref{eq::rv26}) and ends the proof of theorem \ref{le micro resultat}. \rule{2mm}{2mm}

\section{Appendix: Laplace's method}\label{laplace}

For the reader's convenience, we reproduce here Laplace's method. It 
investigates asymptotics as $r$ goes to infinity of integrals involving 
expressions of the form $e^{-rJ}$, $J$ denoting some given function. 
Our interest is to find uniform estimates as 
$x=r(\cos \theta_x,\sin \theta_x)$ lies in a given angle sector $[\theta_1,\theta_2]$.
The proof develops ideas that can be found for a simpler case in \cite{evans}, chapter 4.5.2 page 204.
Lemma \ref{l6.1} below is only used in Step 3 and 4 of the proof of Lemma \ref{supersolution arc}.

\begin{lm}
\label{l6.1}
\textbf{(Uniform asymptotics in a sector $[\theta_1,\theta_2]$) }\\
Define for any $x=r(\cos\theta_x,\sin\theta_x) \in \RR^2$ with $\theta_x \in[0,2\pi)$
$$
F(x)= \lambda_1 \,  e^{\frac{br}{2} \cos(\theta_1-\theta_x)} 
+ \lambda_2 \,  e^{\frac{br}{2} \cos(\theta_2-\theta_x)} 
+ \int_{\theta_1}^{\theta_2} e^{\frac{br}{2} \cos(\theta-\theta_x)} 
f(\theta) \frac{\dd\theta}{2\pi} 
$$
where 
\begin{description}
\item[\tt i.] $b=c_0\cos\alpha >0$, $\lambda_i \in \RR$, 
$\theta_i \in [0,2\pi]$ for $i=1,2$ and $\theta_1< \theta_2$
\item[\tt ii.] $f \in C^1([0,2\pi],\CC)$ is $2\pi$-periodic.
\end{description}
As $r$ goes to infinity, we have the following asymptotics 
uniform in the angular sector $\theta_x \in [\theta_1,\theta_2]$
$$ F(x)=  \lambda_1 \,  e^{\frac{br}{2} \cos(\theta_1-\theta_x)} 
+ \lambda_2 \,  e^{\frac{br}{2} \cos(\theta_2-\theta_x)} 
+e^{\frac{br}{2}} \left(  \frac{f(\theta_x)}{\sqrt{br}} N_0 (x)
+ \frac{R(x)}{r}  \right)$$
where 
$$\displaystyle N_0(x)=\int_{\sqrt{r}g(\theta_1-\theta_x)}^{\sqrt{r}g(\theta_2-\theta_x)} 
e^{-\frac{u^2}{4}} \frac{\dd u}{2\pi} \quad \in [0,1/ \sqrt{\pi}]$$
and 
\begin {equation}
\label{def de g}
g(\theta)=\left\{\begin{array}{ll}
\mbox{ sign}(\theta) \sqrt{2b(1-\cos \theta)} & \quad \mbox{for}\quad \theta\in [-\pi,\pi]\\
\\
\mbox{ sign}(\theta) 2 \sqrt{b}  & \quad \mbox{for}\quad \theta\in \RR\backslash [-\pi,\pi]
\end{array}\right.
\end{equation}
Moreover, there exists a constant $C>0$ such that 
for any $x \in \RR^2$, if $r>1$ and $\theta_x \in [\theta_1,\theta_2]$ then
$|R(x)| \leq C$.
\end{lm}

\noindent\textbf{Proof of Lemma \ref{l6.1}.}\\ 
It is straightforward to check that $g$ defined by \eqref{def de g} is an odd $C^3$-diffeomorphism 
from $[-\pi,\pi]$ to $[-2\sqrt{b},2\sqrt{b}]$ satisfying $g(0)=0$,
$g'(0)=\sqrt{b}$ and $g"(0)=0$.  We have also chosen  to  extend $g$ to the real line by continuity.
However, when we speak about $g^{-1}$,
it has to be understood as the inverse of $g$ on $[-\pi,\pi]$.

Afterwards, for any $x \in \RR^2$, we define
$$
I(x):= \int_{\theta_1}^{\theta_2} e^{\frac{br}{2} \cos(\theta-\theta_x)} f(\theta)\frac{ \dd \theta}{2\pi}=\int_{\theta_1-\theta_x}^{\theta_2-\theta_x} e^{\frac{br}{2} \cos\theta} 
f(\theta+\theta_x) \frac{\dd\theta}{2\pi} 
$$
Assume $\theta_x \in [\theta_1,\theta_2]$. 
In order to get a bound on $I$ uniform in the angle $\theta_x$, we fix some $\delta >0$ and set
$$
\theta_*=\left\{ \begin{array}{ll}
\theta_1-\theta_x & \mbox{ if } \theta_1-\theta_x \geq -\pi+\delta \\
 -\pi + \delta          & \mbox{ otherwise}
\end{array}\right.
\quad 
\theta^*=\left\{ \begin{array}{ll}
\theta_2-\theta_x & \mbox{ if } \theta_2-\theta_x \leq \pi-\delta \\
\pi - \delta          & \mbox{ otherwise}
\end{array}\right.
$$
We then cut the integral $I$ into three parts, integrating between 
$\theta_1-\theta_x$ and $\theta_*$, between $\theta_*$ and $\theta^*$ and 
finally between $\theta^*$ and $\theta_2-\theta_x$. 
We call those three integrals $I_1$, $I_2$ and $I_3$ respectively.

Regarding $I_1$ and $I_3$, $\cos \theta$ can be bounded  in both cases by $\cos (\pi-\delta)$
and $f$ by its $L^{\infty}$ norm on the compact set $[0,2\pi]$.
Thus, there exists a constant $C>0$ such that for any $x \in \RR^2$ 
with  $\theta_x \in [\theta_1,\theta_2]$,
$$I_1+I_3 \leq C e^{\frac{br}{2} \cos(\pi-\delta)}$$
For sufficiently small $\delta>0$, the right hand term decreases exponentially fast 
 and the contribution of $I_1$ and $I_3$ in $I$ is exponentially small 
as $r$ goes to infinity uniformly in $\theta_x$.

Using the change of variables $u=\sqrt{r}g(\theta)$, we rewrite $I_2$ as
$$I_2(x)=\frac{e^{\frac{br}{2}}}{\sqrt{br}} 
\int_{\sqrt{r}g(\theta_*)}^{\sqrt{r}g(\theta^*)}   e^{-\frac{u^2}{4}}
\, h\left(\frac{u}{\sqrt{r}}\right) \frac{\dd u}{2\pi} 
$$
where $h(t)=f(\theta_x+g^{-1}(t))/\sqrt{1-(t^2/(4b))}$ .
Since $h(t)=h(0)+\int_0^t h'(s) \dd s$, we have
$$I_2(x)=e^{\frac{br}{2}} \left( \frac{f(\theta_x)}{\sqrt{br}} N_0^*(x)+ \frac{R(x)}{r} \right)$$
where $N_0^*$ is defined as in lemma \ref{l6.1} with $\theta_*$ or $\theta^*$ when needed, 
but it only changes the desired asymptotics with an exponentially small correction as above.
The remainder term $R$ is defined as 
$$R(x)=r^{\frac{1}{2}}
\int_{\sqrt{r}g(\theta_*)}^{\sqrt{r}g(\theta^*)} e^{-\frac{u^2}{4}} 
\int_0^{\frac{u}{\sqrt{r}}} \frac{h'(s)}{\sqrt{b}}  \dd s \, \frac{\dd u}{2\pi}$$
Since $h$ is smooth, $h'$ is uniformly bounded on $[g(\theta_*),g(\theta^*)]$ and 
the bound only depends on $\delta$. A straight calculation then shows that there exists $C>0$ such that
$$R(x) \leq C \int_{\sqrt{r}g(\theta_*)}^{\sqrt{r}g(\theta^*)} |u| e^{-\frac{u^2}{4}} \dd u 
\leq C  \int_{-\infty}^{+\infty} |u| e^{-\frac{u^2}{4}} \dd u $$
Putting finally $I_1$, $I_2$ and $I_3$ together, we get the desired asymptotics.
\rule{2mm}{2mm}

\bigskip


\end{document}